\documentclass[12pt,letterpaper]{amsart}

\usepackage{amssymb,amsmath,amsthm,amsfonts} 

\usepackage{float}

\usepackage{fullpage}

\usepackage[small]{caption}
\usepackage[lofdepth,lotdepth,captionskip=10pt]{subfig}

\usepackage{pstricks}
\usepackage{epsfig}
\usepackage{pst-node}
\usepackage[dvips={-h tir_____.pfb}]{auto-pst-pdf}

\usepackage[linktocpage=true]{hyperref} 

\usepackage[alphabetic,msc-links,nobysame,lite]{amsrefs} 

\newcommand{\bbC}{\mathbb{C}}

\newcommand{\calC}{\mathcal{C}}

\newcommand{\calP}{\mathcal{P}}

\theoremstyle{plain}

\newtheorem*{lemma*}{Lemma}

\newtheorem*{claim*}{Claim}

\newtheorem*{conjecture*}{Conjecture}

\newtheorem*{corollary*}{Corollary}

\newtheorem*{fact*}{Fact}

\newtheorem*{facts*}{Facts}

\newtheorem*{observation*}{Observation}

\newtheorem*{proposition*}{Proposition}

\newtheorem*{question*}{Question}

\newtheorem*{theorem*}{Theorem}

\newtheorem*{mainindex*}{Main Index Theorem (weak form)}

\theoremstyle{definition}

\newtheorem*{definition*}{Definition}

\newtheorem*{example*}{Example}

\newtheorem*{remark*}{Remark}

\newtheorem*{remarks*}{Remarks}

\theoremstyle{definition}
\newtheorem{exercise}{Problem}

\theoremstyle{plain}
\newtheorem*{exq}{Existence Question}
\newtheorem*{exth}{Existence Theorem}
\newtheorem*{maint}{Uniqueness Theorem}
\newtheorem*{unq}{Uniqueness Question}
\newtheorem*{unt}{Uniqueness Theorem}

\newtheorem{excondition}{Existence Condition}

\begin{document}

\title{Patterns formed by coins}
\author{Andrey M. Mishchenko}
\email{misHchea@umich.edu}
\date{\today}
\thanks{The author was partially supported by NSF grants DMS-0456940, DMS-0555750, DMS-0801029, DMS-1101373.}

\begin{abstract}
This article is a gentle introduction to the mathematical area known as \emph{circle packing}, the study of the kinds of patterns that can be formed by configurations of non-overlapping circles.  The first half of the article is an exposition of the two most important facts about circle packings, (1) that essentially whatever pattern we ask for, we may always arrange circles in that pattern, and (2) that under simple conditions on the pattern, there is an essentially \emph{unique} arrangement of circles in that pattern.  In the second half of the article, we consider related questions, but where we allow the circles to overlap.  The article is written with the idea that no mathematical background should be required to read it.
\end{abstract}

\maketitle

\section{What is a circle packing?}
\label{sec:what is}

\newcommand{\bfit}[1]{\textbf{\emph{#1}}}

Suppose we arrange some coins flat on a table.  We get a picture that looks like this:%
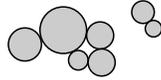
\begin{figure}[H]
\centering
\scalebox{1} 
{
\begin{pspicture}(0,-0.51)(2.06,0.51)
\definecolor{color5245b}{rgb}{0.8,0.8,0.8}
\pscircle[linewidth=0.02,dimen=outer,fillstyle=solid,fillcolor=color5245b](0.23,-0.08){0.23}
\pscircle[linewidth=0.02,dimen=outer,fillstyle=solid,fillcolor=color5245b](0.74,0.11){0.32}
\pscircle[linewidth=0.02,dimen=outer,fillstyle=solid,fillcolor=color5245b](0.94,-0.29){0.14}
\pscircle[linewidth=0.02,dimen=outer,fillstyle=solid,fillcolor=color5245b](1.23,0.04){0.19}
\pscircle[linewidth=0.02,dimen=outer,fillstyle=solid,fillcolor=color5245b](1.25,-0.32){0.19}
\pscircle[linewidth=0.02,dimen=outer,fillstyle=solid,fillcolor=color5245b](1.8,0.35){0.16}
\pscircle[linewidth=0.02,dimen=outer,fillstyle=solid,fillcolor=color5245b](1.94,0.13){0.12}
\end{pspicture} 
}
\end{figure}%
\noindent The coins may be of different sizes, and they may touch, but for now we insist that they do not overlap.  Configurations such as these are traditionally known as \emph{circle packings}:

\begin{definition*}
A \bfit{circle packing} is a collection of disks which don't overlap.
\end{definition*}

We begin by exploring the patterns that can be formed by the disks in a circle packing:

\begin{figure}[H]
\centering
\subfloat
[]
{
\scalebox{1} 
{
\begin{pspicture}(0,-0.62)(1.64,0.62)
\definecolor{color5245b}{rgb}{0.8,0.8,0.8}
\pscircle[linewidth=0.02,dimen=outer,fillstyle=solid,fillcolor=color5245b](0.36,0.26){0.36}
\pscircle[linewidth=0.02,dimen=outer,fillstyle=solid,fillcolor=color5245b](1.03,0.27){0.33}
\pscircle[linewidth=0.02,dimen=outer,fillstyle=solid,fillcolor=color5245b](0.54,-0.34){0.28}
\pscircle[linewidth=0.02,dimen=outer,fillstyle=solid,fillcolor=color5245b](1.01,-0.25){0.21}
\pscircle[linewidth=0.02,dimen=outer,fillstyle=solid,fillcolor=color5245b](1.48,0.42){0.16}
\end{pspicture} 
}
}\qquad\qquad
\subfloat
[]
{
\scalebox{1} 
{
\begin{pspicture}(0,-0.6)(1.6,0.6)
\pscircle[linewidth=0.02,dimen=outer,fillstyle=hlines,hatchwidth=0.01,hatchangle=20.0](0.3,0.32){0.28}
\pscircle[linewidth=0.02,dimen=outer,fillstyle=hlines,hatchwidth=0.01,hatchangle=20.0](0.89,0.25){0.33}
\pscircle[linewidth=0.02,dimen=outer,fillstyle=hlines,hatchwidth=0.01,hatchangle=20.0](0.28,-0.22){0.28}
\pscircle[linewidth=0.02,dimen=outer,fillstyle=hlines,hatchwidth=0.01,hatchangle=20.0](0.81,-0.33){0.27}
\pscircle[linewidth=0.02,dimen=outer,fillstyle=hlines,hatchwidth=0.01,hatchangle=20.0](1.39,0.11){0.21}
\end{pspicture} 
}
}
\caption[Two circle packings which have ``the same pattern'']
{
{\bf Two circle packings which have ``the same pattern,''} in some sense.\label{fig:samesies}
}
\end{figure}
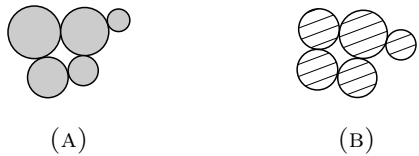

\noindent To make our notion ``pattern'' more precise, we need to know what a \emph{graph} is:

\begin{definition*}
A \bfit{graph} is a collection of \bfit{vertices}, which we usually draw as dots, and \bfit{edges}, which are connections between the dots.
\end{definition*}

\noindent The easiest way to represent a graph is via a drawing.  However, we stress that a graph is an abstract object which is independent of any drawn representation of it.

\begin{figure}[H]
\centering
\subfloat[\label{g1a}]
{
\scalebox{1} 
{
\begin{pspicture}(0,-1.69)(5.04,1.69)
\definecolor{color6218b}{rgb}{0.8,0.8,0.8}
\psline[linewidth=0.02cm,fillcolor=color6218b,dotsize=0.07055555cm 2.0]{*-*}(0.74,0.43)(2.34,-1.17)
\psline[linewidth=0.02cm,fillcolor=color6218b,dotsize=0.07055555cm 2.0]{*-*}(2.34,-1.17)(3.14,0.03)
\psline[linewidth=0.02cm,fillcolor=color6218b,dotsize=0.07055555cm 2.0]{*-*}(3.14,0.03)(0.74,0.43)
\psline[linewidth=0.02cm,fillcolor=color6218b,dotsize=0.07055555cm 2.0]{*-*}(0.74,0.43)(2.54,1.23)
\psline[linewidth=0.02cm,fillcolor=color6218b,dotsize=0.07055555cm 2.0]{*-*}(2.54,1.23)(4.14,0.63)
\psline[linewidth=0.02cm,fillcolor=color6218b,dotsize=0.07055555cm 2.0]{*-*}(4.14,0.63)(3.74,-0.97)
\psline[linewidth=0.02cm,fillcolor=color6218b,dotsize=0.07055555cm 2.0]{*-*}(2.34,-1.17)(3.74,-0.97)
\psline[linewidth=0.02cm,fillcolor=color6218b,dotsize=0.07055555cm 2.0]{*-*}(0.74,0.43)(0.94,-0.97)
\usefont{T1}{ptm}{m}{n}
\rput(0.87,-1.205){$v_1$}
\usefont{T1}{ptm}{m}{n}
\rput(0.45,0.695){$v_2$}
\usefont{T1}{ptm}{m}{n}
\rput(2.53,1.495){$v_3$}
\usefont{T1}{ptm}{m}{n}
\rput(4.55,0.655){$v_4$}
\usefont{T1}{ptm}{m}{n}
\rput(3.97,-1.225){$v_5$}
\usefont{T1}{ptm}{m}{n}
\rput(2.31,-1.465){$v_6$}
\usefont{T1}{ptm}{m}{n}
\rput(3.13,0.255){$v_7$}
\end{pspicture} 
}
}
\qquad
\subfloat[\label{g1b}]
{
\scalebox{1} 
{
\begin{pspicture}(0,-1.68)(5.4826837,1.66)
\definecolor{color6218b}{rgb}{0.8,0.8,0.8}
\usefont{T1}{ptm}{m}{n}
\rput(1.0126834,1.425){$v_1$}
\usefont{T1}{ptm}{m}{n}
\rput(1.3526834,0.225){$v_2$}
\usefont{T1}{ptm}{m}{n}
\rput(2.7726834,1.305){$v_3$}
\usefont{T1}{ptm}{m}{n}
\rput(4.9926834,1.465){$v_4$}
\usefont{T1}{ptm}{m}{n}
\rput(4.8126836,-1.215){$v_5$}
\usefont{T1}{ptm}{m}{n}
\rput(3.1526835,-1.455){$v_6$}
\usefont{T1}{ptm}{m}{n}
\rput(1.1726835,-1.215){$v_7$}
\psbezier[linewidth=0.02,fillcolor=color6218b,dotsize=0.07055555cm 2.0]{*-*}(4.9826837,1.24)(5.3826833,0.8422051)(5.2226834,0.44441018)(5.2226834,0.046615276)(5.2226834,-0.35117963)(5.3826833,-0.76)(4.5826836,-0.947872)
\psbezier[linewidth=0.02,fillcolor=color6218b,dotsize=0.07055555cm 2.0]{*-*}(4.9826837,1.24)(4.1826835,1.24)(4.3826833,0.84)(3.7826834,0.84)(3.1826835,0.84)(3.3826835,1.24)(2.7826834,1.04)
\psbezier[linewidth=0.02,fillcolor=color6218b,dotsize=0.07055555cm 2.0]{*-*}(4.5826836,-0.96)(4.3826833,-0.56)(4.5826836,-0.36)(3.9826834,-0.56)(3.3826835,-0.76)(3.3826835,-0.76)(3.1826835,-1.16)
\psbezier[linewidth=0.02,fillcolor=color6218b,dotsize=0.07055555cm 2.0]{*-*}(3.1826835,-1.16)(2.3826835,-1.56)(2.3826835,-0.96)(1.9826834,-0.96)(1.5826834,-0.96)(1.5826834,-0.36)(0.9826834,-0.96)
\psbezier[linewidth=0.02,fillcolor=color6218b,dotsize=0.07055555cm 2.0]{*-*}(0.9826834,-0.96)(0.58268344,-1.56)(0.0,-1.0838796)(0.38268343,-0.16)(0.76536685,0.76387954)(0.7546162,1.0006288)(1.5826834,0.44)
\psbezier[linewidth=0.02,fillcolor=color6218b,dotsize=0.07055555cm 2.0]{*-*}(1.5826834,0.44)(1.7826835,-0.16)(1.8755766,-0.2528932)(2.1826835,-0.16)(2.4897902,-0.06710678)(2.9826834,-0.16)(3.1826835,-1.16)
\psbezier[linewidth=0.02,fillcolor=color6218b,dotsize=0.07055555cm 2.0]{*-*}(1.5826834,0.44)(2.1826835,0.84)(1.9826834,0.44)(2.3826835,0.64)(2.7826834,0.84)(2.1826835,1.04)(2.7826834,1.04)
\psbezier[linewidth=0.02,fillcolor=color6218b,dotsize=0.07055555cm 2.0]{*-*}(1.5826834,0.44)(1.5826834,1.04)(1.7826835,0.84)(1.7826835,1.24)(1.7826835,1.64)(1.5826834,1.24)(1.3826834,1.44)
\end{pspicture} 
}
}
\caption[Two different representations of the same graph]
{
{\bf Two different representations of the same graph.}  The graph we have drawn here has $7$ vertices.  Here the vertices have been labeled $v_1,v_2,\ldots,v_7$.  Then there is an edge between $v_i$ and $v_j$ in \subref{g1a} whenever there is an edge between $v_i$ and $v_j$ in \subref{g1b} and vice versa.  For example, we do not consider that there is an edge between $v_3$ and $v_5$ in this graph, because we cannot get from $v_3$ to $v_5$ without passing through other vertices.
}
\end{figure}
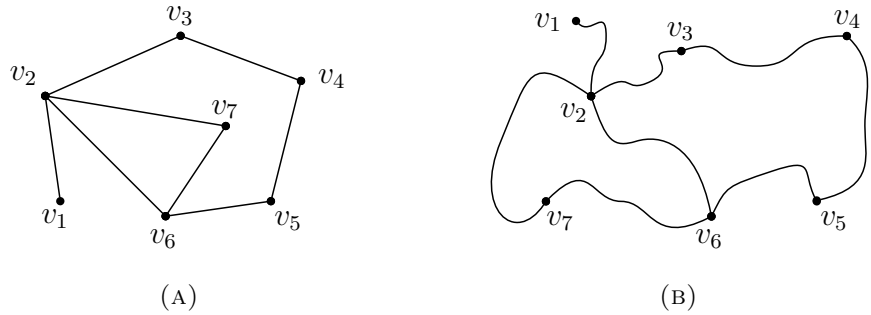

Graphs are all around us.  For example, the \emph{handshake graph} is the graph having a vertex for every person, and having an edge between two people if they have ever shaken hands.  You may be more used to hearing graphs referred to as \emph{networks}.  For example, the \emph{Facebook friend graph} is the graph having a vertex for every Facebook account, and having an edge whenever the two accounts are friended on Facebook.  Of course, the graphs in these two examples both have huge numbers of vertices, and it would be very hard to draw representations of them.  Try the following easier task, to build your intuition:

\begin{exercise}
\label{x:english}
Draw a representation of the following graph:
\begin{itemize}
\item The vertices are the English words for the whole numbers from 1 to 6: one, two, three, four, five, six.
\item Two vertices have an edge connecting them whenever the two corresponding words share a letter.
\end{itemize}
\end{exercise}

The concept of a graph helps us formalize what we mean when we discuss the ``pattern'' formed by the disks of a circle packing.  In particular, we may use a graph to capture the information of ``how the disks of the packing meet'':

\begin{definition*}
Suppose we have a circle packing $\calP$ consisting of the disks $D_1,D_2,D_3,\ldots$.  Then the \bfit{contact graph} of the packing $\calP$ is the graph described as follows:
\begin{itemize}
\item The vertices are the disks $D_1,D_2,D_3,\ldots$ of the packing.
\item We connect two vertices with an edge if, and only if, the corresponding disks touch.
\end{itemize}
\end{definition*}

\noindent We explore some more examples:\vspace{10pt}

\begin{figure}[H]
\centering
\scalebox{1} 
{
\begin{pspicture}(0,-1.14)(3.44,1.14)
\definecolor{color6218b}{rgb}{0.8,0.8,0.8}
\pscircle[linewidth=0.02,dimen=outer,fillstyle=solid,fillcolor=color6218b](0.74,0.4){0.74}
\pscircle[linewidth=0.02,dimen=outer,fillstyle=solid,fillcolor=color6218b](1.87,0.53){0.41}
\pscircle[linewidth=0.02,dimen=outer,fillstyle=solid,fillcolor=color6218b](2.35,-0.03){0.35}
\pscircle[linewidth=0.02,dimen=outer,fillstyle=solid,fillcolor=color6218b](1.43,-0.63){0.51}
\pscircle[linewidth=0.02,dimen=outer,fillstyle=solid,fillcolor=color6218b](2.3,-0.74){0.38}
\pscircle[linewidth=0.02,dimen=outer,fillstyle=solid,fillcolor=color6218b](3.04,0.24){0.4}
\psline[linewidth=0.02cm,fillcolor=color6218b,linestyle=dashed,dash=0.1cm 0.05cm,dotsize=0.07055555cm 2.0]{*-*}(0.74,0.4)(1.87,0.53)
\psline[linewidth=0.02cm,fillcolor=color6218b,linestyle=dashed,dash=0.1cm 0.05cm,dotsize=0.07055555cm 2.0]{*-*}(1.87,0.53)(2.35,-0.03)
\psline[linewidth=0.02cm,fillcolor=color6218b,linestyle=dashed,dash=0.1cm 0.05cm,dotsize=0.07055555cm 2.0]{*-*}(2.35,-0.03)(2.3,-0.74)
\psline[linewidth=0.02cm,fillcolor=color6218b,linestyle=dashed,dash=0.1cm 0.05cm,dotsize=0.07055555cm 2.0]{*-*}(2.3,-0.74)(1.43,-0.63)
\psline[linewidth=0.02cm,fillcolor=color6218b,linestyle=dashed,dash=0.1cm 0.05cm,dotsize=0.07055555cm 2.0]{*-*}(1.43,-0.63)(0.74,0.4)
\psline[linewidth=0.02cm,fillcolor=color6218b,linestyle=dashed,dash=0.1cm 0.05cm,dotsize=0.07055555cm 2.0]{*-*}(2.35,-0.03)(3.04,0.24)
\end{pspicture} 
}\vspace{10pt}
\caption[A circle packing and its contact graph]
{\label{fig:planar rep}
{\bf A circle packing and its contact graph.}  The edges of the graph are drawn as dashed lines.  This is the easiest way to visualize the contact graph of a circle packing, and usually the easiest way to draw a representation of it.
}
\end{figure}
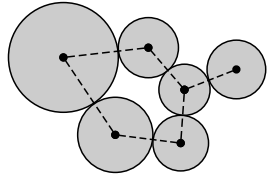

\begin{exercise}
\label{ex333}
Draw the contact graph of the following circle packing:\vspace{10pt}\\ 
\centerline{%
\scalebox{1} 
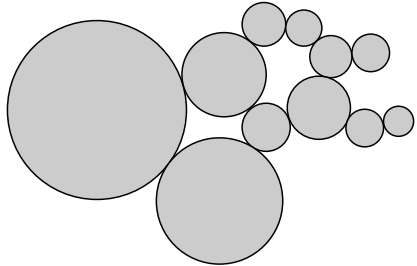
{
\begin{pspicture}(0,-1.75)(4.42,1.75)
\definecolor{color6218b}{rgb}{0.8,0.8,0.8}
\pscircle[linewidth=0.02,dimen=outer,fillstyle=solid,fillcolor=color6218b](1.2,0.31){1.2}
\pscircle[linewidth=0.02,dimen=outer,fillstyle=solid,fillcolor=color6218b](2.89,0.78){0.57}
\pscircle[linewidth=0.02,dimen=outer,fillstyle=solid,fillcolor=color6218b](3.45,0.08){0.33}
\pscircle[linewidth=0.02,dimen=outer,fillstyle=solid,fillcolor=color6218b](2.83,-0.9){0.85}
\pscircle[linewidth=0.02,dimen=outer,fillstyle=solid,fillcolor=color6218b](3.42,1.45){0.3}
\pscircle[linewidth=0.02,dimen=outer,fillstyle=solid,fillcolor=color6218b](3.95,1.4){0.25}
\pscircle[linewidth=0.02,dimen=outer,fillstyle=solid,fillcolor=color6218b](4.31,1.02){0.29}
\pscircle[linewidth=0.02,dimen=outer,fillstyle=solid,fillcolor=color6218b](4.15,0.34){0.43}
\pscircle[linewidth=0.02,dimen=outer,fillstyle=solid,fillcolor=color6218b](4.84,1.07){0.26}
\pscircle[linewidth=0.02,dimen=outer,fillstyle=solid,fillcolor=color6218b](4.76,0.07){0.26}
\pscircle[linewidth=0.02,dimen=outer,fillstyle=solid,fillcolor=color6218b](5.21,0.16){0.21}
\end{pspicture} 
}
}
\end{exercise}

If you did Problem \ref{ex333}, you probably drew the contact graph directly on top of the packing, as we did in Figure \ref{fig:planar rep}, but it is not always necessary to do this:

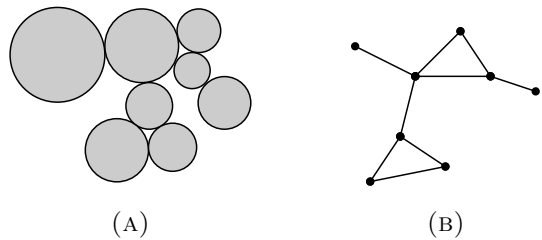
\begin{figure}[H]
\centering
\subfloat[\label{amm1}]
{
\scalebox{1} 
{
\begin{pspicture}(0,-1.17)(3.22,1.17)
\definecolor{color6218b}{rgb}{0.8,0.8,0.8}
\pscircle[linewidth=0.02,dimen=outer,fillstyle=solid,fillcolor=color6218b](0.64,0.53){0.64}
\pscircle[linewidth=0.02,dimen=outer,fillstyle=solid,fillcolor=color6218b](1.76,0.65){0.5}
\pscircle[linewidth=0.02,dimen=outer,fillstyle=solid,fillcolor=color6218b](1.86,-0.15){0.32}
\pscircle[linewidth=0.02,dimen=outer,fillstyle=solid,fillcolor=color6218b](2.43,0.32){0.25}
\pscircle[linewidth=0.02,dimen=outer,fillstyle=solid,fillcolor=color6218b](2.52,0.85){0.3}
\pscircle[linewidth=0.02,dimen=outer,fillstyle=solid,fillcolor=color6218b](1.43,-0.74){0.43}
\pscircle[linewidth=0.02,dimen=outer,fillstyle=solid,fillcolor=color6218b](2.17,-0.7){0.33}
\pscircle[linewidth=0.02,dimen=outer,fillstyle=solid,fillcolor=color6218b](2.86,-0.11){0.36}
\end{pspicture} 
}
}
\qquad
\subfloat[\label{amm2}]
{
\scalebox{1} 
{
\begin{pspicture}(0,-1.01)(2.41,1.01)
\definecolor{color6218b}{rgb}{0.8,0.8,0.8}
\psline[linewidth=0.02cm,fillcolor=color6218b,dotsize=0.07055555cm 2.0]{*-*}(0.2,-1.0)(0.6,-0.4)
\psline[linewidth=0.02cm,fillcolor=color6218b,dotsize=0.07055555cm 2.0]{*-*}(0.6,-0.4)(1.2,-0.8)
\psline[linewidth=0.02cm,fillcolor=color6218b,dotsize=0.07055555cm 2.0]{*-*}(1.2,-0.8)(0.2,-1.0)
\psline[linewidth=0.02cm,fillcolor=color6218b,dotsize=0.07055555cm 2.0]{*-*}(0.6,-0.4)(0.8,0.4)
\psline[linewidth=0.02cm,fillcolor=color6218b,dotsize=0.07055555cm 2.0]{*-*}(0.8,0.4)(1.4,1.0)
\psline[linewidth=0.02cm,fillcolor=color6218b,dotsize=0.07055555cm 2.0]{*-*}(1.4,1.0)(1.8,0.4)
\psline[linewidth=0.02cm,fillcolor=color6218b,dotsize=0.07055555cm 2.0]{*-*}(1.8,0.4)(0.8,0.4)
\psline[linewidth=0.02cm,fillcolor=color6218b,dotsize=0.07055555cm 2.0]{*-*}(1.8,0.4)(2.4,0.2)
\psline[linewidth=0.02cm,fillcolor=color6218b,dotsize=0.07055555cm 2.0]{*-*}(0.8,0.4)(0.0,0.8)
\end{pspicture} 
}
}
\caption[Another circle packing and its contact graph]
{
{\bf Another circle packing and its contact graph.}  Here we have drawn the contact graph separately from the packing, to again emphasize that the graph is an abstract object distinct from any drawn representation of it.
}
\end{figure}

\begin{exercise}
Verify that the graph shown in Figure \ref{amm2} is the contact graph of the circle packing shown in Figure \ref{amm1}.
\end{exercise}

Experience leads us to conclude:

\begin{observation*}
If we start with a circle packing, it is not so hard to draw its contact graph.
\end{observation*}

\noindent In other words, it is easy to ``go from'' a circle packing to its contact graph.  Then we might ask, when can we go ``in the other direction?''  This is our segue into the next section.

\section{Existence of circle packings}
\label{lay:ex}

To build intuition, we start this section with an exercise:

\begin{exercise}
\label{x:existence}
Draw circle packings having the following contact graphs:\\ \\
\centerline{
\scalebox{1} 
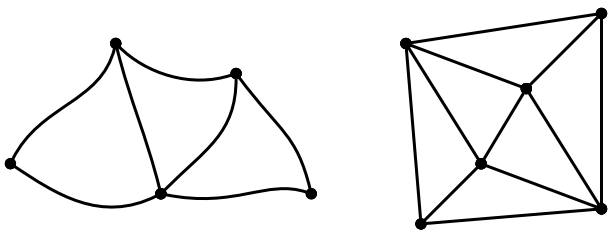
{
\begin{pspicture}(0,-1.22)(4.02,1.22)
\psbezier[linewidth=0.04,dotsize=0.07055555cm 2.0]{*-*}(1.4,1.2)(1.6,0.4)(1.8,0.0)(2.0,-0.8)
\psbezier[linewidth=0.04,dotsize=0.07055555cm 2.0]{*-*}(2.0,-0.8)(1.2,-1.2)(0.6,-0.8)(0.0,-0.4)
\psbezier[linewidth=0.04,dotsize=0.07055555cm 2.0]{*-*}(0.0,-0.4)(0.4,0.4)(1.2,0.4)(1.4,1.2)
\psbezier[linewidth=0.04,dotsize=0.07055555cm 2.0]{*-*}(1.4,1.2)(1.8,0.8)(2.4,0.6)(3.0,0.8)
\psbezier[linewidth=0.04,dotsize=0.07055555cm 2.0]{*-*}(3.0,0.8)(3.0,0.0)(2.6,-0.2)(2.0,-0.8)
\psbezier[linewidth=0.04,dotsize=0.07055555cm 2.0]{*-*}(2.0,-0.8)(3.0,-1.0)(3.4,-0.6)(4.0,-0.8)
\psbezier[linewidth=0.04,dotsize=0.07055555cm 2.0]{*-*}(4.0,-0.8)(3.8,0.0)(3.6,0.0)(3.0,0.8)
\end{pspicture} 
}
\qquad
\scalebox{1} 
{
\begin{pspicture}(0,-1.42)(2.62,1.42)
\psline[linewidth=0.04cm,dotsize=0.07055555cm 2.0]{*-*}(0.0,1.0)(0.2,-1.4)
\psline[linewidth=0.04cm,dotsize=0.07055555cm 2.0]{*-*}(0.2,-1.4)(2.6,-1.2)
\psline[linewidth=0.04cm,dotsize=0.07055555cm 2.0]{*-*}(2.6,-1.2)(2.6,1.4)
\psline[linewidth=0.04cm,dotsize=0.07055555cm 2.0]{*-*}(2.6,1.4)(0.0,1.0)
\psline[linewidth=0.04cm,dotsize=0.07055555cm 2.0]{*-*}(0.0,1.0)(1.6,0.4)
\psline[linewidth=0.04cm,dotsize=0.07055555cm 2.0]{*-*}(1.6,0.4)(2.6,1.4)
\psline[linewidth=0.04cm,dotsize=0.07055555cm 2.0]{*-*}(1.6,0.4)(2.6,-1.2)
\psline[linewidth=0.04cm,dotsize=0.07055555cm 2.0]{*-*}(2.6,-1.2)(1.0,-0.6)
\psline[linewidth=0.04cm,dotsize=0.07055555cm 2.0]{*-*}(1.0,-0.6)(0.2,-1.4)
\psline[linewidth=0.04cm,dotsize=0.07055555cm 2.0]{*-*}(1.0,-0.6)(1.6,0.4)
\psline[linewidth=0.04cm,dotsize=0.07055555cm 2.0]{*-*}(1.0,-0.6)(0.0,1.0)
\end{pspicture} 
}
}
\end{exercise}

\noindent This leads us to consider the following natural question:

\begin{exq}
If we start with a graph $G$, is there always a circle packing having $G$ as its contact graph?
\end{exq}

\noindent For simplicity, we will consider only \emph{finite} graphs, meaning graphs which have finitely many vertices and edges.  We first observe that circle packings' contact graphs never have \emph{loops}, and never have \emph{repeated edges}:%
\begin{figure}[H]
\centering
\scalebox{1} 
{
\begin{pspicture}(0,-2.0994906)(4.42,2.0994904)
\psbezier[linewidth=0.04,dotsize=0.07055555cm 2.0]{*-*}(0.82,-0.08069181)(1.22,-0.08069181)(2.02,-0.2806918)(2.22,0.5193082)
\psbezier[linewidth=0.04,dotsize=0.07055555cm 2.0]{*-*}(0.82,-0.08069181)(0.82,-0.8806918)(1.22,-0.8806918)(1.62,-0.8806918)
\psbezier[linewidth=0.04,dotsize=0.07055555cm 2.0]{*-*}(2.22,0.5193082)(2.82,0.3193082)(3.5012093,0.595914)(3.22,-0.08069181)(2.9387908,-0.75729764)(2.42,-0.08069181)(2.22,0.5193082)
\psbezier[linewidth=0.04,dotsize=0.07055555cm 2.0]{*-*}(2.22,0.5193082)(1.62,-0.08069181)(1.0730537,-0.018576821)(1.42,0.9193082)(1.7669462,1.8571932)(1.6045878,1.3075136)(2.22,0.5193082)
\psbezier[linewidth=0.04,dotsize=0.07055555cm 2.0]{*-*}(0.82,-0.08069181)(0.82,-0.8806918)(0.83291256,1.759126)(1.82,1.9193082)(2.8070874,2.0794904)(3.82,0.5193082)(3.82,-0.08069181)(3.82,-0.6806918)(1.82,-0.48069182)(1.62,-0.8806918)(1.42,-1.2806919)(2.0380435,-1.4915843)(2.22,-1.4806918)(2.4019563,-1.4697993)(2.42,-1.0806918)(2.82,-1.2806919)(3.22,-1.4806918)(3.570997,-2.0794904)(3.62,-1.0806918)
\usefont{T1}{ptm}{m}{n}
\rput(2.51,0.7843082){$v_1$}
\usefont{T1}{ptm}{m}{n}
\rput(0.45,-0.035691805){$v_2$}
\usefont{T1}{ptm}{m}{n}
\rput(1.25,-1.1356918){$v_3$}
\usefont{T1}{ptm}{m}{n}
\rput(3.93,-0.8956918){$v_4$}
\end{pspicture} 
}
\caption[A graph which \emph{does} have loops and repeated edges]
{
{\bf A graph which \emph{does} have loops and repeated edges.}  The vertices $v_2$ and $v_3$ have two edges connecting them, an example of what is meant by a \emph{repeated edge}.  The vertex $v_1$ has two loops on it: a \emph{loop} is an edge from a vertex to itself.
}
\end{figure}
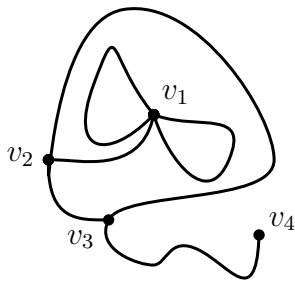%
\noindent A graph without loops and repeated edges is called \emph{simple}.  We observe:

\begin{excondition}
\label{cond rice}
If $G$ is the contact graph of a circle packing, then $G$ is simple.
\end{excondition}

Next, we make a more subtle observation:

\begin{fact*}
Not every graph can be drawn without its edges crossing each other.
\end{fact*}

\noindent The idea of \emph{edges crossing} is illustrated in the following example:
\begin{figure}[H]
\centering
\subfloat[\label{gg1}]
{
\scalebox{1} 
{
\begin{pspicture}(0,-1.52)(3.22,1.52)
\psbezier[linewidth=0.04,dotsize=0.15055555cm 2.0]{*-*}(0.0,1.5)(0.0,0.7)(0.4,-1.3)(0.4,-1.3)
\psbezier[linewidth=0.04,dotsize=0.15055555cm 2.0]{*-*}(0.4,-1.3)(1.2,-1.5)(1.6,-1.1)(2.4,-1.1)
\psbezier[linewidth=0.04,dotsize=0.15055555cm 2.0]{*-*}(2.4,-1.1)(2.6,-0.3)(3.2,0.7)(2.8,1.5)
\psbezier[linewidth=0.04,dotsize=0.15055555cm 2.0]{*-*}(2.8,1.5)(1.8,1.5)(0.8,1.1)(0.0,1.5)
\psbezier[linewidth=0.04,dotsize=0.15055555cm 2.0]{*-*}(0.4,-1.3)(0.8,-0.5)(2.8,0.7)(2.8,1.5)
\psbezier[linewidth=0.04,dotsize=0.15055555cm 2.0]{*-*}(0.0,1.5)(0.8,0.1)(2.4,-1.1)(2.4,-1.1)
\end{pspicture} 
}
}\qquad\qquad
\subfloat[\label{gg2}]
{
\scalebox{1} 
{
\begin{pspicture}(0,-1.62)(3.82,1.62)
\psbezier[linewidth=0.04,dotsize=0.15055555cm 2.0]{*-*}(0.2,0.6)(0.0,-0.2)(0.2,-1.4)(0.2,-1.4)
\psbezier[linewidth=0.04,dotsize=0.15055555cm 2.0]{*-*}(0.2,-1.4)(1.0,-1.6)(1.4,-1.2)(2.2,-1.2)
\psbezier[linewidth=0.04,dotsize=0.15055555cm 2.0]{*-*}(2.2,-1.2)(2.4,-0.4)(1.8,-0.4)(1.4,0.4)
\psbezier[linewidth=0.04,dotsize=0.15055555cm 2.0]{*-*}(1.4,0.4)(1.0,0.6)(0.8,0.8)(0.2,0.6)
\psbezier[linewidth=0.04,dotsize=0.15055555cm 2.0]{*-*}(0.2,-1.4)(0.4,-0.6)(1.4,-0.4)(1.4,0.4)
\psbezier[linewidth=0.04,dotsize=0.15055555cm 2.0]{*-*}(0.2,0.6)(1.2,1.6)(3.8,1.2)(2.2,-1.2)
\end{pspicture} 
}
}
\caption[Two representations of the same graph, one planar and one not]
{
{\bf Two representations of the same graph, one planar and one not.}  The graph represented here has four vertices.  In \subref{gg1}, two edges of our graph cross.  This does not happen in \subref{gg2}.  The drawing in \subref{gg2} is called a \emph{planar representation} for this graph.
}
\end{figure}
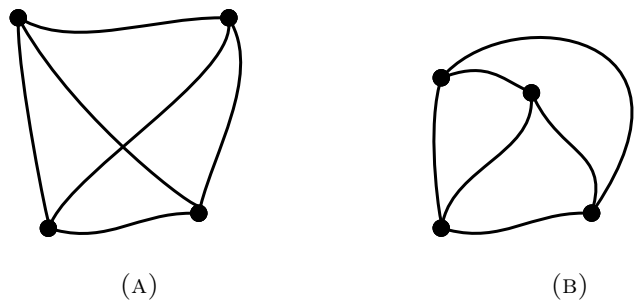

\noindent A graph is called \emph{planar} if it can be drawn (in the plane) without edge crossings.  (You can think of ``the plane'' as an infinite flat piece of paper, that goes on forever in all directions.)  To convince yourself that not all graphs are planar, try the following exercise:

\begin{exercise}
\label{ex:k5}
Try to draw a representation of the following graph, so that the edges don't cross each other:
\begin{itemize}
\item It has five vertices, and
\item \emph{every pair} of vertices has an edge connecting them.
\end{itemize}
In other words, draw five dots on a piece of paper, and try to draw a connection between every pair of the dots so that these connections don't cross.  (It is impossible, but try to convince yourself of this, rather than taking it on faith.)
\end{exercise}

On the other hand, suppose that the graph $G$ is the contact graph of some circle packing.  Then we may draw a representation of $G$ as we did in Figure \ref{fig:planar rep}, by putting the vertices corresponding to our disks at the centers of the disks, and connecting the vertices of touching disks by straight line segments.  We conclude:

\begin{excondition}
\label{cond rice 2}
If $G$ is the contact graph of a circle packing, then $G$ is planar.
\end{excondition}

\noindent Thus for example, there is no circle packing whose graph is the one described in Exercise \ref{ex:k5}.\medskip

Before moving on, do the following quick exercise.  We will use your answer soon.

\begin{exercise}
\label{newplanar}
Draw a simple planar graph different from any of the ones we have seen so far.
\end{exercise}

An amazing fact is that Existence Conditions \ref{cond rice} and \ref{cond rice 2} are the \emph{only conditions} needed to guarantee the existence of a circle packing in whatever pattern you like:

\begin{exth}
If $G$ is a simple planar graph, then there is some circle packing having $G$ as its contact graph.
\end{exth}

\noindent Thus, without seeing your answer to Problem \ref{newplanar}, I know that it is possible to do the following:

\begin{exercise}
Draw a circle packing whose contact graph is the one you drew in Problem \ref{newplanar}.
\end{exercise}

The Existence Question for circle packings was first asked by Paul Koebe, and he answered it himself in 1936, by proving the above existence theorem.  The proof uses \emph{complex analysis}, the theory of calculus on the complex numbers $\bbC = \{ a + bi \}$.  It turns out that the field of circle packing is a kind of discrete analog of complex analysis.  This is a major motivating factor for the study of circle packings.  The reference for Koebe's original paper is \cite{koebe-1936}.

\section{Rigidity of circle packings}
\label{lay:un}

Different circle packings can have the same contact graph, for example, the packings in Figure \ref{fig:samesies} on p.\ \pageref{fig:samesies}.  In this section, we explore the question of whether there is a contact graph so that \emph{only one} circle packing has that contact graph.

Of course, if we start with a circle packing and, imaging that our disks are glued together wherever they touch, pick all of the disks up, move them around, and then set them back down, we will get a new circle packing having the same contact graph.  This leads us to explore the notion of \emph{similar packings}.

Recall from elementary school geometry that two triangles are called \emph{similar} if one is a copy of the first, but scaled up or down and moved around:

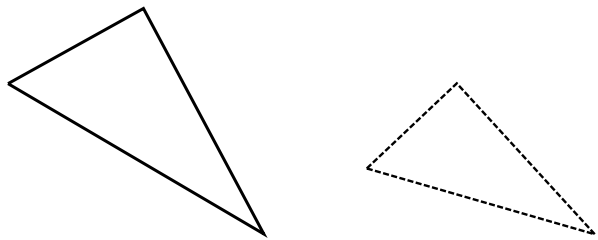
\begin{figure}[H]
\centering
\subfloat
{
\scalebox{1} 
{
\begin{pspicture}(0,-1.52)(3.42,1.52)
\definecolor{color711b}{rgb}{0.8,0.8,0.8}
\psline[linewidth=0.04,fillcolor=color711b](0.0,0.5)(1.8,1.5)(3.4,-1.5)(0.0,0.5)
\end{pspicture} 
}
}
\qquad
\subfloat
{
\scalebox{.8} 
{
\begin{pspicture}(0,-1.2749563)(3.8095777,1.2749565)
\definecolor{color711b}{rgb}{0.8,0.8,0.8}
\psline[linewidth=0.04,linestyle=dashed,dash=0.1cm 0.04cm,fillcolor=color711b](0.0,-0.15992126)(1.4960351,1.2549565)(3.7895775,-1.2549565)(0.0,-0.15992126)
\end{pspicture} 
}
}
\caption[Two similar triangles]
{
\label{simtri}
{\bf Two similar triangles.}
}
\end{figure}

\noindent We can define two packings to be similar in an analogous way:

\begin{definition*}
Two circle packings are called \bfit{similar} if it is possible to get from one to the other by combining the following four operations:
\begin{itemize}
\item sliding,
\item rotating,
\item reflecting (``flipping''), and
\item scaling (``stretching'').
\end{itemize}
\end{definition*}

\noindent Remember that the operations are always applied to a \emph{whole packing at once}, not to individual disks one at a time!  We consider similar packings to be ``essentially the same'':

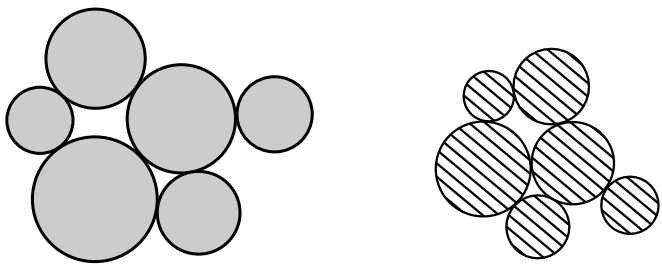
\begin{figure}[H]
\centering
\subfloat
{
\scalebox{1} 
{
\begin{pspicture}(0,-1.7)(4.1,1.7)
\definecolor{color627b}{rgb}{0.8,0.8,0.8}
\pscircle[linewidth=0.04,dimen=outer,fillstyle=solid,fillcolor=color627b](1.19,-0.85){0.85}
\pscircle[linewidth=0.04,dimen=outer,fillstyle=solid,fillcolor=color627b](2.34,0.22){0.74}
\pscircle[linewidth=0.04,dimen=outer,fillstyle=solid,fillcolor=color627b](2.57,-1.03){0.57}
\pscircle[linewidth=0.04,dimen=outer,fillstyle=solid,fillcolor=color627b](3.58,0.28){0.52}
\pscircle[linewidth=0.04,dimen=outer,fillstyle=solid,fillcolor=color627b](0.46,0.2){0.46}
\pscircle[linewidth=0.04,dimen=outer,fillstyle=solid,fillcolor=color627b](1.2,1.02){0.68}
\end{pspicture} 
}
}
\qquad
\subfloat
{
\scalebox{.76} 
{
\begin{pspicture}(0,-1.9428114)(4.344149,1.9116791)
\rput{-39.2085}(0.45435393,0.69008803){\pscircle[linewidth=0.04,dimen=outer,fillstyle=hlines,hatchwidth=0.04,hatchangle=0.0](1.1959457,-0.29279333){0.85}}
\rput{-39.2085}(0.7427119,1.70395){\pscircle[linewidth=0.04,dimen=outer,fillstyle=hlines,hatchwidth=0.04,hatchangle=0.0](2.7634184,-0.19066894){0.74}}
\rput{-39.2085}(1.3091098,1.0662936){\pscircle[linewidth=0.04,dimen=outer,fillstyle=hlines,hatchwidth=0.04,hatchangle=0.0](2.1514537,-1.3046255){0.57}}
\rput{-39.2085}(1.4337012,2.1692822){\pscircle[linewidth=0.04,dimen=outer,fillstyle=hlines,hatchwidth=0.04,hatchangle=0.0](3.762162,-0.92803675){0.52}}
\rput{-39.2085}(-0.3295775,1.0391865){\pscircle[linewidth=0.04,dimen=outer,fillstyle=hlines,hatchwidth=0.04,hatchangle=0.0](1.2940562,0.98226523){0.46}}
\rput{-39.2085}(-0.18971288,1.7670615){\pscircle[linewidth=0.04,dimen=outer,fillstyle=hlines,hatchwidth=0.04,hatchangle=0.0](2.385804,1.1498561){0.68}}
\end{pspicture} 
}
}
\caption[Two circle packings which are ``essentially the same'']
{
\label{simcp}
{\bf Two circle packings which are ``essentially the same.''}
}
\end{figure}

\noindent The reason for this is that similar packings \emph{always} have the same contact graph.

On the other hand, there may still be some circle packing $\calP$, so that \emph{every} packing which has the same contact graph as $\calP$ turns out to be similar to $\calP$.  In other words:

\begin{unq}
When is there an \emph{essentially} unique circle packing having some given contact graph?
\end{unq}

\noindent For every circle packing we have seen so far, it is possible to find an essentially different packing, meaning a packing which is \emph{not} similar to the original one, which has the same contact graph.  Convince yourself:

\begin{exercise}
Look back, and pick one of the circle packings we have seen so far, call it $\calP$.  Draw a new circle packing which has the same contact graph as your chosen packing $\calP$, but which is not similar to $\calP$.
\end{exercise}

In fact, it is not too hard to believe the following:

\begin{fact*}
If a circle packing $\calP$ is made up of finitely many disks, then there are circle packings having the same contact graph as $\calP$, but which are not similar to $\calP$.
\end{fact*}

\noindent So, we turn our attention to \emph{infinite packings}, packings which have infinitely many disks: 

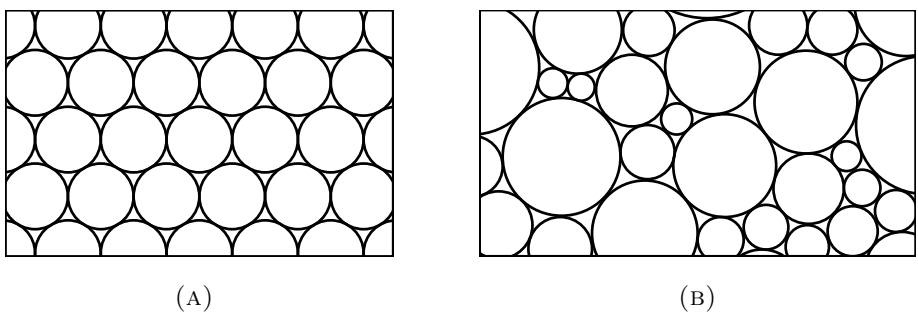
\begin{figure}[H]
\centering
{
\subfloat[\label{infaaa1}\label{inf1}]
{
\scalebox{0.436797}{
\begin{pspicture*}(0.1,-0.1)(11.9,7.4)
\psframe[dimen=middle,linewidth=0.0824](0.1,-0.1)(11.9,7.4)
\pscircle[dimen=middle,linewidth=0.0824](00,0){1}
\pscircle[dimen=middle,linewidth=0.0824](02,0){1}
\pscircle[dimen=middle,linewidth=0.0824](04,0){1}
\pscircle[dimen=middle,linewidth=0.0824](06,0){1}
\pscircle[dimen=middle,linewidth=0.0824](08,0){1}
\pscircle[dimen=middle,linewidth=0.0824](10,0){1}
\pscircle[dimen=middle,linewidth=0.0824](12,0){1}%

\pscircle[dimen=middle,linewidth=0.0824](01,1.73){1}
\pscircle[dimen=middle,linewidth=0.0824](03,1.73){1}
\pscircle[dimen=middle,linewidth=0.0824](05,1.73){1}
\pscircle[dimen=middle,linewidth=0.0824](07,1.73){1}
\pscircle[dimen=middle,linewidth=0.0824](09,1.73){1}
\pscircle[dimen=middle,linewidth=0.0824](11,1.73){1}
\pscircle[dimen=middle,linewidth=0.0824](13,1.73){1}%

\pscircle[dimen=middle,linewidth=0.0824](00,3.46){1}
\pscircle[dimen=middle,linewidth=0.0824](02,3.46){1}
\pscircle[dimen=middle,linewidth=0.0824](04,3.46){1}
\pscircle[dimen=middle,linewidth=0.0824](06,3.46){1}
\pscircle[dimen=middle,linewidth=0.0824](08,3.46){1}
\pscircle[dimen=middle,linewidth=0.0824](10,3.46){1}
\pscircle[dimen=middle,linewidth=0.0824](12,3.46){1}%

\pscircle[dimen=middle,linewidth=0.0824](01,5.19){1}
\pscircle[dimen=middle,linewidth=0.0824](03,5.19){1}
\pscircle[dimen=middle,linewidth=0.0824](05,5.19){1}
\pscircle[dimen=middle,linewidth=0.0824](07,5.19){1}
\pscircle[dimen=middle,linewidth=0.0824](09,5.19){1}
\pscircle[dimen=middle,linewidth=0.0824](11,5.19){1}
\pscircle[dimen=middle,linewidth=0.0824](13,5.19){1}%

\pscircle[dimen=middle,linewidth=0.0824](00,6.928){1}
\pscircle[dimen=middle,linewidth=0.0824](02,6.928){1}
\pscircle[dimen=middle,linewidth=0.0824](04,6.928){1}
\pscircle[dimen=middle,linewidth=0.0824](06,6.928){1}
\pscircle[dimen=middle,linewidth=0.0824](08,6.928){1}
\pscircle[dimen=middle,linewidth=0.0824](10,6.928){1}
\pscircle[dimen=middle,linewidth=0.0824](12,6.928){1}%
\end{pspicture*}
}
}
\qquad
\subfloat[\label{infaaa2}\label{inf2}]
{
\scalebox{.9} 
{
\begin{pspicture*}(7.56,1.12)(1.12,-2.52)
\psframe[linewidth=0.04,dimen=middle](7.56,1.12)(1.12,-2.52)
\pscircle[linewidth=0.04,dimen=outer](1.01,0.27){1.01}
\pscircle[linewidth=0.04,dimen=outer](2.57,0.83){0.67}
\pscircle[linewidth=0.04,dimen=outer](3.37,-0.07){0.55}
\pscircle[linewidth=0.04,dimen=outer](2.2,0.04){0.24}
\pscircle[linewidth=0.04,dimen=outer](2.62,-0.02){0.22}
\pscircle[linewidth=0.04,dimen=outer](2.33,-1.05){0.89}
\pscircle[linewidth=0.04,dimen=outer](3.62,0.82){0.4}
\pscircle[linewidth=0.04,dimen=outer](4.56,0.28){0.72}
\pscircle[linewidth=0.04,dimen=outer](3.6,-0.98){0.42}
\pscircle[linewidth=0.04,dimen=outer](4.03,-0.49){0.25}
\pscircle[linewidth=0.04,dimen=outer](1.01,-1.19){0.47}
\pscircle[linewidth=0.04,dimen=outer](4.74,-1.18){0.78}
\pscircle[linewidth=0.04,dimen=outer](3.56,-2.18){0.8}
\pscircle[linewidth=0.04,dimen=outer](1.4,-2.06){0.52}
\pscircle[linewidth=0.04,dimen=outer](2.31,-2.41){0.49}
\pscircle[linewidth=0.04,dimen=outer](4.68,-2.28){0.36}
\pscircle[linewidth=0.04,dimen=outer](5.35,-2.09){0.35}
\pscircle[linewidth=0.04,dimen=outer](5.98,-1.52){0.54}
\pscircle[linewidth=0.04,dimen=outer](5.94,-0.24){0.78}
\pscircle[linewidth=0.04,dimen=outer](6.33,0.83){0.39}
\pscircle[linewidth=0.04,dimen=outer](6.79,0.35){0.29}
\pscircle[linewidth=0.04,dimen=outer](7.7,-0.58){1.04}
\pscircle[linewidth=0.04,dimen=outer](6.54,-1.04){0.24}
\pscircle[linewidth=0.04,dimen=outer](6.77,-1.51){0.29}
\pscircle[linewidth=0.04,dimen=outer](6.63,-2.15){0.39}
\pscircle[linewidth=0.04,dimen=outer](5.96,-2.38){0.34}
\pscircle[linewidth=0.04,dimen=outer](7.27,-1.85){0.33}
\pscircle[linewidth=0.04,dimen=outer](7.42,1.2){0.78}
\pscircle[linewidth=0.04,dimen=outer](5.53,0.89){0.45}
\pscircle[linewidth=0.04,dimen=outer](7.36,-2.7){0.56}
\pscircle[linewidth=0.04,dimen=outer](5.29,-3.01){0.61}
\pscircle[linewidth=0.04,dimen=outer](4.48,2.3){1.32}
\end{pspicture*} 
}
}
}
\caption[Pieces of infinite circle packings]
{\label{fig:infaaa} \label{fig:inf}
{\bf Pieces of infinite circle packings.}  Of course, it is impossible to draw the entirety of a packing that goes on forever.  However, sometimes, it is still possible to completely describe the packing.  For example, the packing shown in \subref{inf1} is the infinite packing by disks which all have the same size, so that each disk has exactly 6 neighbors.  This is known as the ``penny packing.''  We leave what happens in \subref{inf2} to the imagination.
}
\end{figure}

\noindent We might be tempted to believe that if a packing ``goes on forever in all directions,'' then it is essentially the only packing having its contact graph.  Unfortunately, this does not turn out to be true.  For example:

\begin{figure}[H]
\centering
\subfloat[\label{sq1}]
{
\scalebox{.6}
{
\begin{pspicture*}(-.1,-.1)(9.1,5.1)
\psframe[dimen=middle](-.1,-.1)(9.1,5.1)
\pscircle[dimen=middle](0,0){.5}
\pscircle[dimen=middle](1,0){.5}
\pscircle[dimen=middle](2,0){.5}
\pscircle[dimen=middle](3,0){.5}
\pscircle[dimen=middle](4,0){.5}
\pscircle[dimen=middle](5,0){.5}
\pscircle[dimen=middle](6,0){.5}
\pscircle[dimen=middle](7,0){.5}
\pscircle[dimen=middle](8,0){.5}
\pscircle[dimen=middle](9,0){.5}%

\pscircle[dimen=middle](0,1){.5}
\pscircle[dimen=middle](1,1){.5}
\pscircle[dimen=middle](2,1){.5}
\pscircle[dimen=middle](3,1){.5}
\pscircle[dimen=middle](4,1){.5}
\pscircle[dimen=middle](5,1){.5}
\pscircle[dimen=middle](6,1){.5}
\pscircle[dimen=middle](7,1){.5}
\pscircle[dimen=middle](8,1){.5}
\pscircle[dimen=middle](9,1){.5}%

\pscircle[dimen=middle](0,2){.5}
\pscircle[dimen=middle](1,2){.5}
\pscircle[dimen=middle](2,2){.5}
\pscircle[dimen=middle](3,2){.5}
\pscircle[dimen=middle](4,2){.5}
\pscircle[dimen=middle](5,2){.5}
\pscircle[dimen=middle](6,2){.5}
\pscircle[dimen=middle](7,2){.5}
\pscircle[dimen=middle](8,2){.5}
\pscircle[dimen=middle](9,2){.5}%

\pscircle[dimen=middle](0,3){.5}
\pscircle[dimen=middle](1,3){.5}
\pscircle[dimen=middle](2,3){.5}
\pscircle[dimen=middle](3,3){.5}
\pscircle[dimen=middle](4,3){.5}
\pscircle[dimen=middle](5,3){.5}
\pscircle[dimen=middle](6,3){.5}
\pscircle[dimen=middle](7,3){.5}
\pscircle[dimen=middle](8,3){.5}
\pscircle[dimen=middle](9,3){.5}%

\pscircle[dimen=middle](0,4){.5}
\pscircle[dimen=middle](1,4){.5}
\pscircle[dimen=middle](2,4){.5}
\pscircle[dimen=middle](3,4){.5}
\pscircle[dimen=middle](4,4){.5}
\pscircle[dimen=middle](5,4){.5}
\pscircle[dimen=middle](6,4){.5}
\pscircle[dimen=middle](7,4){.5}
\pscircle[dimen=middle](8,4){.5}
\pscircle[dimen=middle](9,4){.5}%

\pscircle[dimen=middle](0,5){.5}
\pscircle[dimen=middle](1,5){.5}
\pscircle[dimen=middle](2,5){.5}
\pscircle[dimen=middle](3,5){.5}
\pscircle[dimen=middle](4,5){.5}
\pscircle[dimen=middle](5,5){.5}
\pscircle[dimen=middle](6,5){.5}
\pscircle[dimen=middle](7,5){.5}
\pscircle[dimen=middle](8,5){.5}
\pscircle[dimen=middle](9,5){.5}
\end{pspicture*}
}
}
\qquad
\subfloat[\label{sq2}]
{
\scalebox{.6}
{
\begin{pspicture*}(-.1,-.1)(9.1,5.1)
\psframe[dimen=middle](-.1,-.1)(9.1,5.1)
\pscircle[dimen=middle](0,0){.5}
\pscircle[dimen=middle](1,0){.5}
\pscircle[dimen=middle](2,0){.5}
\pscircle[dimen=middle](3,0){.5}
\pscircle[dimen=middle](4,0){.5}
\pscircle[dimen=middle](4.978,.2){.5}
\pscircle[dimen=middle](5.978,.2){.5}
\pscircle[dimen=middle](6.978,.2){.5}
\pscircle[dimen=middle](7.978,.2){.5}
\pscircle[dimen=middle](8.978,.2){.5}%

\pscircle[dimen=middle](0,1){.5}
\pscircle[dimen=middle](1,1){.5}
\pscircle[dimen=middle](2,1){.5}
\pscircle[dimen=middle](3,1){.5}
\pscircle[dimen=middle](4,1){.5}
\pscircle[dimen=middle](4.978,1.2){.5}
\pscircle[dimen=middle](5.978,1.2){.5}
\pscircle[dimen=middle](6.978,1.2){.5}
\pscircle[dimen=middle](7.978,1.2){.5}
\pscircle[dimen=middle](8.978,1.2){.5}%

\pscircle[dimen=middle](0,2){.5}
\pscircle[dimen=middle](1,2){.5}
\pscircle[dimen=middle](2,2){.5}
\pscircle[dimen=middle](3,2){.5}
\pscircle[dimen=middle](4,2){.5}
\pscircle[dimen=middle](4.978,2.2){.5}
\pscircle[dimen=middle](5.978,2.2){.5}
\pscircle[dimen=middle](6.978,2.2){.5}
\pscircle[dimen=middle](7.978,2.2){.5}
\pscircle[dimen=middle](8.978,2.2){.5}%

\pscircle[dimen=middle](0,3){.5}
\pscircle[dimen=middle](1,3){.5}
\pscircle[dimen=middle](2,3){.5}
\pscircle[dimen=middle](3,3){.5}
\pscircle[dimen=middle](4,3){.5}
\pscircle[dimen=middle](4.978,3.2){.5}
\pscircle[dimen=middle](5.978,3.2){.5}
\pscircle[dimen=middle](6.978,3.2){.5}
\pscircle[dimen=middle](7.978,3.2){.5}
\pscircle[dimen=middle](8.978,3.2){.5}%

\pscircle[dimen=middle](0,4){.5}
\pscircle[dimen=middle](1,4){.5}
\pscircle[dimen=middle](2,4){.5}
\pscircle[dimen=middle](3,4){.5}
\pscircle[dimen=middle](4,4){.5}
\pscircle[dimen=middle](4.978,4.2){.5}
\pscircle[dimen=middle](5.978,4.2){.5}
\pscircle[dimen=middle](6.978,4.2){.5}
\pscircle[dimen=middle](7.978,4.2){.5}
\pscircle[dimen=middle](8.978,4.2){.5}%

\pscircle[dimen=middle](0,5){.5}
\pscircle[dimen=middle](1,5){.5}
\pscircle[dimen=middle](2,5){.5}
\pscircle[dimen=middle](3,5){.5}
\pscircle[dimen=middle](4,5){.5}
\pscircle[dimen=middle](4.978,5.2){.5}
\pscircle[dimen=middle](5.978,5.2){.5}
\pscircle[dimen=middle](6.978,5.2){.5}
\pscircle[dimen=middle](7.978,5.2){.5}
\pscircle[dimen=middle](8.978,5.2){.5}
\end{pspicture*}
}
}
\caption[Two infinite packings which are not similar, but which have the same contact graph]
{
\label{fig:sq}
{\bf Two infinite packings which are not similar, but which have the same contact graph.}
}
\end{figure}
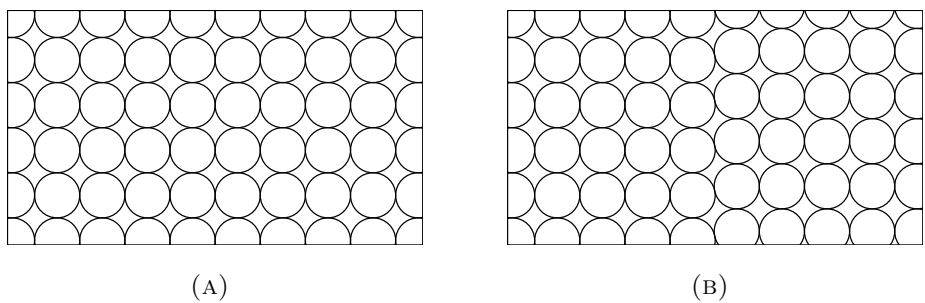

\begin{exercise}
Describe the contact graph of the packings in Figure \ref{fig:sq}.  Can you find a third packing which has the same contact graph, but which is not similar to either of \subref{sq1} and \subref{sq2}?
\end{exercise}

In Figure \ref{fig:sq}, we ``perturbed'' some disks of the packing, leaving others where they were, to get a new, fundamentally different packing having the same contact graph.  On the other hand, it does not seem so easy to perturb in a similar way the disks of either packing in Figure \ref{fig:infaaa}.  This is because of a special property of the contact graph:

\begin{definition*}
Suppose that $\calP$ is an infinite circle packing.  Represent its contact graph $G$ by putting a vertex at every circle's center, and connecting vertices of touching circles by straight segments.  Suppose that then the graph $G$ ``cuts the plane into triangles.''  In that case we say that the packing $\calP$ is \bfit{triangulated}.
\end{definition*} 

\begin{figure}[H]
\centering
\scalebox{0.48533}{
\begin{pspicture*}(0.1,-0.3)(11.9,7.4)
\psframe[dimen=middle,linewidth=0.0824](0.1,-0.3)(11.9,7.4)%

\psdots[dotsize=0.5](00,0)
\psdots[dotsize=0.5](02,0)
\psdots[dotsize=0.5](04,0)
\psdots[dotsize=0.5](06,0)
\psdots[dotsize=0.5](08,0)
\psdots[dotsize=0.5](10,0)
\psdots[dotsize=0.5](12,0)%

\psdots[dotsize=0.5](01,1.73)
\psdots[dotsize=0.5](03,1.73)
\psdots[dotsize=0.5](05,1.73)
\psdots[dotsize=0.5](07,1.73)
\psdots[dotsize=0.5](09,1.73)
\psdots[dotsize=0.5](11,1.73)
\psdots[dotsize=0.5](13,1.73)%

\psdots[dotsize=0.5](00,3.46)
\psdots[dotsize=0.5](02,3.46)
\psdots[dotsize=0.5](04,3.46)
\psdots[dotsize=0.5](06,3.46)
\psdots[dotsize=0.5](08,3.46)
\psdots[dotsize=0.5](10,3.46)
\psdots[dotsize=0.5](12,3.46)%

\psdots[dotsize=0.5](01,5.19)
\psdots[dotsize=0.5](03,5.19)
\psdots[dotsize=0.5](05,5.19)
\psdots[dotsize=0.5](07,5.19)
\psdots[dotsize=0.5](09,5.19)
\psdots[dotsize=0.5](11,5.19)
\psdots[dotsize=0.5](13,5.19)%

\psdots[dotsize=0.5](00,6.928)
\psdots[dotsize=0.5](02,6.928)
\psdots[dotsize=0.5](04,6.928)
\psdots[dotsize=0.5](06,6.928)
\psdots[dotsize=0.5](08,6.928)
\psdots[dotsize=0.5](10,6.928)
\psdots[dotsize=0.5](12,6.928)%

\pscircle[dimen=middle,linewidth=0.0324](00,0){1}
\pscircle[dimen=middle,linewidth=0.0324](02,0){1}
\pscircle[dimen=middle,linewidth=0.0324](04,0){1}
\pscircle[dimen=middle,linewidth=0.0324](06,0){1}
\pscircle[dimen=middle,linewidth=0.0324](08,0){1}
\pscircle[dimen=middle,linewidth=0.0324](10,0){1}
\pscircle[dimen=middle,linewidth=0.0324](12,0){1}%

\pscircle[dimen=middle,linewidth=0.0324](01,1.73){1}
\pscircle[dimen=middle,linewidth=0.0324](03,1.73){1}
\pscircle[dimen=middle,linewidth=0.0324](05,1.73){1}
\pscircle[dimen=middle,linewidth=0.0324](07,1.73){1}
\pscircle[dimen=middle,linewidth=0.0324](09,1.73){1}
\pscircle[dimen=middle,linewidth=0.0324](11,1.73){1}
\pscircle[dimen=middle,linewidth=0.0324](13,1.73){1}%

\pscircle[dimen=middle,linewidth=0.0324](00,3.46){1}
\pscircle[dimen=middle,linewidth=0.0324](02,3.46){1}
\pscircle[dimen=middle,linewidth=0.0324](04,3.46){1}
\pscircle[dimen=middle,linewidth=0.0324](06,3.46){1}
\pscircle[dimen=middle,linewidth=0.0324](08,3.46){1}
\pscircle[dimen=middle,linewidth=0.0324](10,3.46){1}
\pscircle[dimen=middle,linewidth=0.0324](12,3.46){1}%

\pscircle[dimen=middle,linewidth=0.0324](01,5.19){1}
\pscircle[dimen=middle,linewidth=0.0324](03,5.19){1}
\pscircle[dimen=middle,linewidth=0.0324](05,5.19){1}
\pscircle[dimen=middle,linewidth=0.0324](07,5.19){1}
\pscircle[dimen=middle,linewidth=0.0324](09,5.19){1}
\pscircle[dimen=middle,linewidth=0.0324](11,5.19){1}
\pscircle[dimen=middle,linewidth=0.0324](13,5.19){1}%

\pscircle[dimen=middle,linewidth=0.0324](00,6.928){1}
\pscircle[dimen=middle,linewidth=0.0324](02,6.928){1}
\pscircle[dimen=middle,linewidth=0.0324](04,6.928){1}
\pscircle[dimen=middle,linewidth=0.0324](06,6.928){1}
\pscircle[dimen=middle,linewidth=0.0324](08,6.928){1}
\pscircle[dimen=middle,linewidth=0.0324](10,6.928){1}
\pscircle[dimen=middle,linewidth=0.0324](12,6.928){1}%

\psline[linewidth=0.09,linestyle=dashed,dash=0.2cm 0.05cm](-1,0)(15,0)
\psline[linewidth=0.09,linestyle=dashed,dash=0.2cm 0.05cm](-1,1.73)(15,1.73)
\psline[linewidth=0.09,linestyle=dashed,dash=0.2cm 0.05cm](-1,3.46)(15,3.46)
\psline[linewidth=0.09,linestyle=dashed,dash=0.2cm 0.05cm](-1,5.19)(15,5.19)
\psline[linewidth=0.09,linestyle=dashed,dash=0.2cm 0.05cm](-1,6.928)(15,6.928)%

\pnode(-5,-1.73){Y}
\pnode(-3,-1.73){Z}
\pnode(-1,-1.73){A}
\pnode(1,-1.73){B}
\pnode(3,-1.73){C}
\pnode(5,-1.73){D}
\pnode(7,-1.73){E}
\pnode(9,-1.73){F}
\pnode(11,-1.73){G}
\pnode(13,-1.73){H}
\pnode(15,-1.73){I}
\pnode(17,-1.73){J}
\pnode(19,-1.73){K}
\SpecialCoor
\psline[linewidth=0.09,linestyle=dashed,dash=0.2cm 0.05cm](Y)([angle=60,nodesep=15]Y)
\psline[linewidth=0.09,linestyle=dashed,dash=0.2cm 0.05cm](Z)([angle=60,nodesep=15]Z)
\psline[linewidth=0.09,linestyle=dashed,dash=0.2cm 0.05cm](A)([angle=60,nodesep=15]A)
\psline[linewidth=0.09,linestyle=dashed,dash=0.2cm 0.05cm](B)([angle=60,nodesep=15]B)
\psline[linewidth=0.09,linestyle=dashed,dash=0.2cm 0.05cm](C)([angle=60,nodesep=15]C)
\psline[linewidth=0.09,linestyle=dashed,dash=0.2cm 0.05cm](D)([angle=60,nodesep=15]D)
\psline[linewidth=0.09,linestyle=dashed,dash=0.2cm 0.05cm](E)([angle=60,nodesep=15]E)
\psline[linewidth=0.09,linestyle=dashed,dash=0.2cm 0.05cm](F)([angle=60,nodesep=15]F)
\psline[linewidth=0.09,linestyle=dashed,dash=0.2cm 0.05cm](G)([angle=60,nodesep=15]G)
\psline[linewidth=0.09,linestyle=dashed,dash=0.2cm 0.05cm](H)([angle=60,nodesep=15]H)
\psline[linewidth=0.09,linestyle=dashed,dash=0.2cm 0.05cm](I)([angle=60,nodesep=15]I)
\psline[linewidth=0.09,linestyle=dashed,dash=0.2cm 0.05cm](J)([angle=60,nodesep=15]J)
\psline[linewidth=0.09,linestyle=dashed,dash=0.2cm 0.05cm](K)([angle=60,nodesep=15]K)%

\psline[linewidth=0.09,linestyle=dashed,dash=0.2cm 0.05cm](Y)([angle=120,nodesep=15]Y)
\psline[linewidth=0.09,linestyle=dashed,dash=0.2cm 0.05cm](Z)([angle=120,nodesep=15]Z)
\psline[linewidth=0.09,linestyle=dashed,dash=0.2cm 0.05cm](A)([angle=120,nodesep=15]A)
\psline[linewidth=0.09,linestyle=dashed,dash=0.2cm 0.05cm](B)([angle=120,nodesep=15]B)
\psline[linewidth=0.09,linestyle=dashed,dash=0.2cm 0.05cm](C)([angle=120,nodesep=15]C)
\psline[linewidth=0.09,linestyle=dashed,dash=0.2cm 0.05cm](D)([angle=120,nodesep=15]D)
\psline[linewidth=0.09,linestyle=dashed,dash=0.2cm 0.05cm](E)([angle=120,nodesep=15]E)
\psline[linewidth=0.09,linestyle=dashed,dash=0.2cm 0.05cm](F)([angle=120,nodesep=15]F)
\psline[linewidth=0.09,linestyle=dashed,dash=0.2cm 0.05cm](G)([angle=120,nodesep=15]G)
\psline[linewidth=0.09,linestyle=dashed,dash=0.2cm 0.05cm](H)([angle=120,nodesep=15]H)
\psline[linewidth=0.09,linestyle=dashed,dash=0.2cm 0.05cm](I)([angle=120,nodesep=15]I)
\psline[linewidth=0.09,linestyle=dashed,dash=0.2cm 0.05cm](J)([angle=120,nodesep=15]J)
\psline[linewidth=0.09,linestyle=dashed,dash=0.2cm 0.05cm](K)([angle=120,nodesep=15]K)
\end{pspicture*}
}
\caption[The contact graph of the penny packing]
{\label{contactpenny}
{\bf The contact graph of the penny packing.}  We see that the contact graph, drawn in the way we described, cuts the plane into triangles.
}
\end{figure}
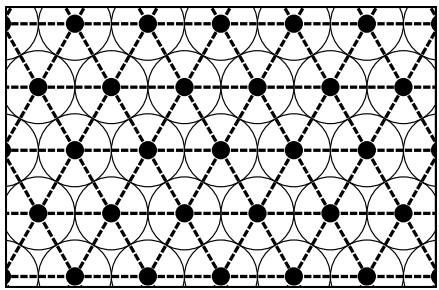

\begin{exercise}
\label{ex:draw cg}
Verify that the packing in Figure \ref{infaaa2} appears to be triangulated as well.  (Of course, we don't really know what happens outside of the rectangular region we have drawn.)
\end{exercise}

\noindent The following amazing fact turns out to be true:

\begin{unt}
If two circle packings have the same contact graph and are are both triangulated, then they are similar.  In other words, if a packing $\calP$ is triangulated and has contact graph $G$, then $\calP$ is essentially the \emph{only} circle packing having contact graph $G$.
\end{unt}

\noindent A alternative way to phrase this theorem is to say that triangulated circle packings are \emph{rigid}, meaning that we cannot modify a triangulated circle packing without changing its contact graph, except by rotating, scaling, etc.  Compare this to the situation in Figure \ref{fig:sq}, where it \emph{is} possible to modify the packings significantly while keeping the contact graph the same.  Furthermore, considering the example of Figure \ref{fig:sq}, the following is not hard to believe:

\begin{fact*}
Suppose a circle packing $\calP$ is not triangulated.  Then there are circle packings having the same contact graph as $\calP$, but which are not similar to $\calP$.
\end{fact*}

\noindent In other words, the condition that $\calP$ is triangulated is \emph{necessary} for $\calP$ to be rigid.  The Uniqueness Theorem says that this condition is also \emph{sufficient}.  Thus the uniqueness guaranteed by the Uniqueness Theorem is the best that we could possibly hope for.

As we mentioned earlier, Koebe initiated the study of circle packings in 1936, and proved the Existence Theorem for circle packings.  He also managed to prove a different kind of uniqueness theorem for circle packings: instead of considering packings in the plane, as we are, he considered packings on the surface of a sphere.  The truth or falsehood of the above uniqueness theorem was unknown for a long time after Koebe, until a proof was finally found in 1991 by Oded Schramm.  The reference for the original article by Schramm is \cite{MR1076089}.

\section{Overlapping disks}
\label{lay:ov}

We next ask ourselves, what happens if we allow the disks to overlap?  It will be helpful to have some vocabulary available to us:

\begin{definition*}
A \bfit{disk configuration} is any collection of disks, so that no disk in the collection is completely contained inside of another.
\end{definition*}%
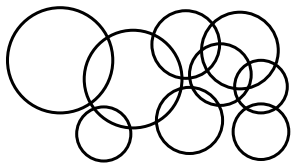
\begin{figure}[H]
\centering
\scalebox{1} 
{
\begin{pspicture}(0,-1.04)(3.78,1.04)
\pscircle[linewidth=0.04,dimen=outer](1.69,0.07){0.67}
\pscircle[linewidth=0.04,dimen=outer](2.39,0.53){0.47}
\pscircle[linewidth=0.04,dimen=outer](2.85,0.13){0.43}
\pscircle[linewidth=0.04,dimen=outer](2.44,-0.48){0.46}
\pscircle[linewidth=0.04,dimen=outer](3.11,0.45){0.53}
\pscircle[linewidth=0.04,dimen=outer](3.39,-0.03){0.37}
\pscircle[linewidth=0.04,dimen=outer](3.39,-0.63){0.39}
\pscircle[linewidth=0.04,dimen=outer](0.72,0.32){0.72}
\pscircle[linewidth=0.04,dimen=outer](1.3,-0.66){0.38}
\end{pspicture} 
}
\caption[Our first example of a disk configuration]
{
{\bf Our first example of a disk configuration.}
}
\end{figure}%
\noindent When discussing disk configurations, it will be helpful to have a notion of \emph{overlap angle}:

\begin{figure}[H]
\centering
\subfloat
{
\scalebox{.9} 
{
\begin{pspicture}(0,-1.39)(3.66,1.4)
\definecolor{color740b}{rgb}{0.8,0.8,0.8}
\pscircle[linewidth=0.04,dimen=outer](1.16,-0.15){1.16}
\pscircle[linewidth=0.04,dimen=outer](2.73,-0.46){0.93}
\psline[linewidth=0.02cm,fillcolor=color740b,arrowsize=0.05291667cm 2.0,arrowlength=1.4,arrowinset=0.4]{->}(2.22,0.31)(3.24,0.91)
\psline[linewidth=0.02cm,fillcolor=color740b,arrowsize=0.05291667cm 2.0,arrowlength=1.4,arrowinset=0.4]{->}(2.2,0.31)(1.8,1.39)
\psarc[linewidth=0.02,fillcolor=color740b](2.21,0.28){0.41}{30.963757}{110.22486}
\usefont{T1}{ptm}{m}{n}
\rput(2.4014552,0.935){$\theta$}
\end{pspicture} 
}
}\quad
\subfloat
{
\scalebox{.9} 
{
\begin{pspicture}(0,-1.46)(4.06,1.468)
\definecolor{color740b}{rgb}{0.8,0.8,0.8}
\pscircle[linewidth=0.04,dimen=outer](1.16,-0.3){1.16}
\pscircle[linewidth=0.04,dimen=outer](3.13,0.05){0.93}
\psline[linewidth=0.02cm,fillcolor=color740b,arrowsize=0.05291667cm 2.0,arrowlength=1.4,arrowinset=0.4]{->}(2.22,0.16)(2.26,1.24)
\psline[linewidth=0.02cm,fillcolor=color740b,arrowsize=0.05291667cm 2.0,arrowlength=1.4,arrowinset=0.4]{->}(2.2,0.16)(1.8,1.24)
\psarc[linewidth=0.02,fillcolor=color740b](2.21,0.13){0.41}{86.98721}{110.22486}
\usefont{T1}{ptm}{m}{n}
\rput(0.7214551,1.145){$\theta$}
\psbezier[linewidth=0.016,fillcolor=color740b,arrowsize=0.05291667cm 2.0,arrowlength=1.4,arrowinset=0.4]{<-}(2.14,0.6)(2.12,1.36)(1.78,1.46)(1.58,1.36)(1.38,1.26)(1.2,1.14)(0.96,1.14)
\end{pspicture} 
}
}\quad
\subfloat
{
\scalebox{.9} 
{
\begin{pspicture}(0,-1.35)(3.2818944,1.36)
\definecolor{color740b}{rgb}{0.8,0.8,0.8}
\pscircle[linewidth=0.04,dimen=outer](1.16,-0.19){1.16}
\pscircle[linewidth=0.04,dimen=outer](1.57,-0.38){0.93}
\psline[linewidth=0.02cm,fillcolor=color740b,arrowsize=0.05291667cm 2.0,arrowlength=1.4,arrowinset=0.4]{->}(2.22,0.27)(2.9,-0.29)
\psline[linewidth=0.02cm,fillcolor=color740b,arrowsize=0.05291667cm 2.0,arrowlength=1.4,arrowinset=0.4]{->}(2.2,0.27)(1.8,1.35)
\psarc[linewidth=0.02,fillcolor=color740b](2.21,0.24){0.41}{324.92624}{110.22486}
\usefont{T1}{ptm}{m}{n}
\rput(2.7614552,0.615){$\theta$}
\end{pspicture} 
}
}
\caption[The overlap angles between several pairs of disks]
{
{\bf The overlap angles between several pairs of disks.}  In each case, we have chosen a point where the boundary circles meet, and drawn the tangent rays to the circles pointing outside of the disks.  Then the \emph{overlap angle} $\theta$ is the angle between these rays.  We note two facts: (1) the overlap angle gets bigger as the disks ``move closer together,'' or gets smaller as the disks ``move apart,'' and (2) the overlap angle is the same at both points of intersection of the boundary circles.
}
\end{figure}
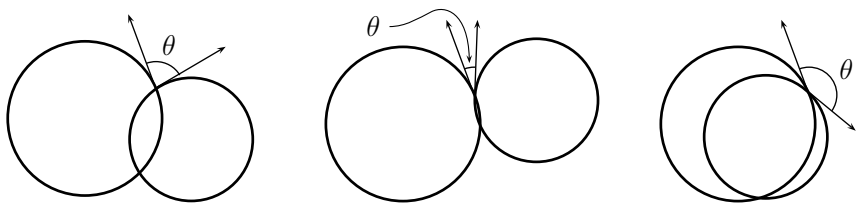

The \bfit{contact graph} of a disk configuration is defined the same way as the \bfit{contact graph} of a circle packing.  In the case of disk configurations, however, we label each edge of the contact graph with the overlap angle between those two disks.  For example:%
\begin{figure}[H]
\centering
\scalebox{1} 
{
\begin{pspicture}(0,-1.9329102)(7.98291,1.9329102)
\definecolor{color966b}{rgb}{0.8,0.8,0.8}
\pscircle[linewidth=0.04,dimen=outer](1.8410156,-0.20548828){1.06}
\pscircle[linewidth=0.04,dimen=outer](3.2610157,0.11451172){1.06}
\pscircle[linewidth=0.04,dimen=outer](5.1010156,0.59451175){0.86}
\pscircle[linewidth=0.04,dimen=outer](2.8210156,-0.9854883){0.72}
\pscircle[linewidth=0.04,dimen=outer](5.4910154,-0.5954883){0.61}
\pscircle[linewidth=0.04,dimen=outer](6.6210155,0.6545117){0.68}
\pscircle[linewidth=0.04,dimen=outer](6.9410157,0.33451173){0.42}
\psline[linewidth=0.04cm,fillcolor=color966b,linestyle=dashed,dash=0.1cm 0.05cm,dotsize=0.07055555cm 2.0]{*-*}(2.9610157,-1.2654883)(1.7610157,-0.26548827)
\psline[linewidth=0.04cm,fillcolor=color966b,linestyle=dashed,dash=0.1cm 0.05cm,dotsize=0.07055555cm 2.0]{*-*}(1.7610157,-0.26548827)(3.3610156,0.13451172)
\psline[linewidth=0.04cm,fillcolor=color966b,linestyle=dashed,dash=0.1cm 0.05cm,dotsize=0.07055555cm 2.0]{*-*}(3.3610156,0.13451172)(2.9610157,-1.2654883)
\psline[linewidth=0.04cm,fillcolor=color966b,linestyle=dashed,dash=0.1cm 0.05cm,dotsize=0.07055555cm 2.0]{*-*}(3.3610156,0.13451172)(5.1610155,0.53451174)
\psline[linewidth=0.04cm,fillcolor=color966b,linestyle=dashed,dash=0.1cm 0.05cm,dotsize=0.07055555cm 2.0]{*-*}(5.1610155,0.53451174)(5.5610156,-0.6654883)
\psline[linewidth=0.04cm,fillcolor=color966b,linestyle=dashed,dash=0.1cm 0.05cm,dotsize=0.07055555cm 2.0]{*-*}(5.1610155,0.53451174)(6.301016,0.59451175)
\psline[linewidth=0.04cm,fillcolor=color966b,linestyle=dashed,dash=0.1cm 0.05cm,dotsize=0.07055555cm 2.0]{*-*}(6.301016,0.59451175)(7.1610155,0.11451172)
\psbezier[linewidth=0.02,fillcolor=color966b,arrowsize=0.05291667cm 2.0,arrowlength=1.4,arrowinset=0.4]{<-}(2.3610156,-0.8454883)(2.1610155,-1.1854882)(2.2810156,-1.5654883)(2.0210156,-1.6254883)(1.7610157,-1.6854882)(1.6410156,-1.5054883)(1.3210156,-1.5454882)
\psbezier[linewidth=0.02,fillcolor=color966b,arrowsize=0.05291667cm 2.0,arrowlength=1.4,arrowinset=0.4]{<-}(2.5810156,-0.005488281)(2.4410157,0.43451172)(2.0810156,1.1145117)(1.8610157,1.1545117)(1.6410156,1.1945118)(1.6410156,0.9545117)(1.1010156,0.85451174)
\psbezier[linewidth=0.02,fillcolor=color966b,arrowsize=0.05291667cm 2.0,arrowlength=1.4,arrowinset=0.4]{<-}(3.1810157,-0.6854883)(3.5410156,-1.0054883)(3.3852792,-1.4363576)(3.5010157,-1.5654883)(3.6167521,-1.6946189)(3.8210156,-1.7654883)(4.0410156,-1.7454883)
\psbezier[linewidth=0.02,fillcolor=color966b,arrowsize=0.05291667cm 2.0,arrowlength=1.4,arrowinset=0.4]{<-}(4.0010157,0.37451172)(4.0010157,0.97451174)(4.2010155,1.0345117)(4.0810156,1.2145118)(3.9610157,1.3945117)(3.7410157,1.3945117)(3.5610156,1.5545117)
\psbezier[linewidth=0.02,fillcolor=color966b,arrowsize=0.05291667cm 2.0,arrowlength=1.4,arrowinset=0.4]{<-}(5.4410157,-0.40548828)(5.1210155,-0.5054883)(4.6410155,-0.4254883)(4.4410157,-0.6054883)(4.2410154,-0.7854883)(4.2010155,-1.0454882)(4.5410156,-1.2254883)
\psbezier[linewidth=0.02,fillcolor=color966b,arrowsize=0.05291667cm 2.0,arrowlength=1.4,arrowinset=0.4]{<-}(5.6210155,0.6345117)(5.5610156,0.9945117)(5.6210155,1.3345118)(5.841016,1.5545117)(6.0610156,1.7745117)(5.5410156,1.6345117)(5.1810155,1.7145118)
\psbezier[linewidth=0.02,fillcolor=color966b,arrowsize=0.05291667cm 2.0,arrowlength=1.4,arrowinset=0.4]{<-}(6.7810154,0.27451172)(6.5410156,0.13451172)(6.3947873,-0.04444163)(6.341016,-0.18548828)(6.287244,-0.32653493)(6.4410157,-0.52548826)(6.7210155,-0.52548826)
\usefont{T1}{ptm}{m}{n}
\rput(0.7224707,0.8195117){$90^\circ$}
\usefont{T1}{ptm}{m}{n}
\rput(1.0224707,-1.5404882){$90^\circ$}
\usefont{T1}{ptm}{m}{n}
\rput(4.3424706,-1.7604883){$90^\circ$}
\usefont{T1}{ptm}{m}{n}
\rput(3.4724708,1.6995118){$0^\circ$}
\usefont{T1}{ptm}{m}{n}
\rput(4.9024707,-1.3204883){$30^\circ$}
\usefont{T1}{ptm}{m}{n}
\rput(4.9724708,1.7395117){$0^\circ$}
\usefont{T1}{ptm}{m}{n}
\rput(7.0924706,-0.5404883){$150^\circ$}
\end{pspicture} 
}
\caption[The contact graph and overlap angles of this disk configuration]
{
{\bf The contact graph and overlap angles of a disk configuration.}\label{jeff}
}
\end{figure}
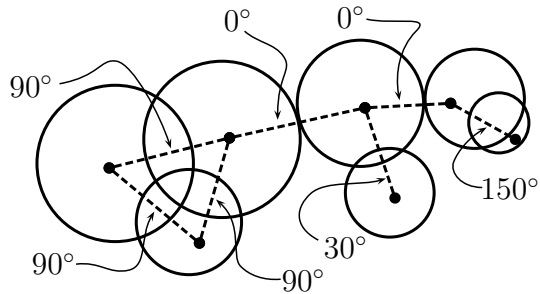

\noindent Try the following problem to get a handle on these new concepts:

\begin{exercise}
Draw a disk configuration having the following contact graph and (approximate) overlap angles:
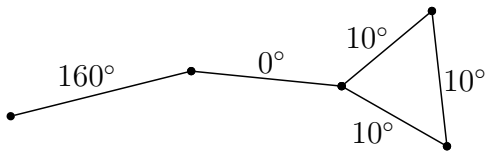
\begin{figure}[H]
\centering
\scalebox{1} 
{
\begin{pspicture}(0,-0.91)(6.8418946,0.91)
\definecolor{color1652b}{rgb}{0.8,0.8,0.8}
\psline[linewidth=0.02cm,fillcolor=color1652b,dotsize=0.07055555cm 2.0]{*-*}(0.0,-0.5)(2.4,0.1)
\psline[linewidth=0.02cm,fillcolor=color1652b,dotsize=0.07055555cm 2.0]{*-*}(2.4,0.1)(4.4,-0.1)
\psline[linewidth=0.02cm,fillcolor=color1652b,dotsize=0.07055555cm 2.0]{*-*}(4.4,-0.1)(5.6,0.9)
\psline[linewidth=0.02cm,fillcolor=color1652b,dotsize=0.07055555cm 2.0]{*-*}(5.6,0.9)(5.8,-0.9)
\psline[linewidth=0.02cm,fillcolor=color1652b,dotsize=0.07055555cm 2.0]{*-*}(5.8,-0.9)(4.4,-0.1)
\usefont{T1}{ptm}{m}{n}
\rput(1.011455,0.045){$160^\circ$}
\usefont{T1}{ptm}{m}{n}
\rput(3.471455,0.225){$0^\circ$}
\usefont{T1}{ptm}{m}{n}
\rput(4.741455,0.545){$10^\circ$}
\usefont{T1}{ptm}{m}{n}
\rput(6.0414553,-0.015){$10^\circ$}
\usefont{T1}{ptm}{m}{n}
\rput(4.821455,-0.715){$10^\circ$}
\end{pspicture} 
}
\end{figure}
\end{exercise}

In our new setting, we consider that the ``pattern'' of a disk configuration consists of both its contact graph, and its overlap angles.  We ask ourselves the following natural question:

\begin{question*}
Given a graph $G$, so that its edges are labeled with some angles, when does there exist a disk configuration having $G$ as its contact graph with those overlap angles?
\end{question*}

\noindent This question turns out to be much harder than the existence question for circle packings.  For our first example, consider the following:

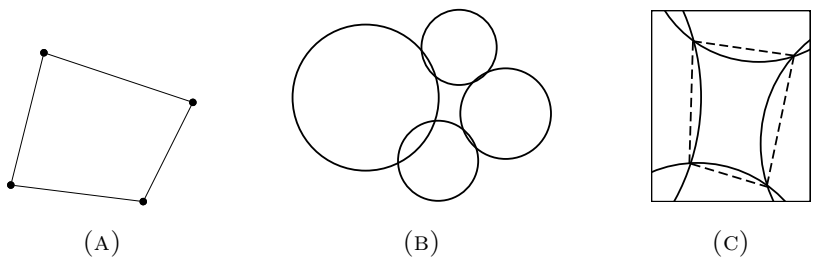
\begin{figure}[H]
\centering
\subfloat[\label{art1}]
{
\scalebox{1.1} 
{
\begin{pspicture}(0,-0.905)(2.205,0.905)
\definecolor{color1768b}{rgb}{0.8,0.8,0.8}
\psline[linewidth=0.01cm,fillcolor=color1768b,dotsize=0.07055555cm 2.0]{*-*}(0.4,0.9)(0.0,-0.7)
\psline[linewidth=0.01cm,fillcolor=color1768b,dotsize=0.07055555cm 2.0]{*-*}(0.0,-0.7)(1.6,-0.9)
\psline[linewidth=0.01cm,fillcolor=color1768b,dotsize=0.07055555cm 2.0]{*-*}(1.6,-0.9)(2.2,0.3)
\psline[linewidth=0.01cm,fillcolor=color1768b,dotsize=0.07055555cm 2.0]{*-*}(2.2,0.3)(0.4,0.9)
\end{pspicture} 
}
}\qquad
\subfloat[\label{art2}]
{
\scalebox{.41} 
{
\begin{pspicture}(0,-3.13)(8.44,3.13)
\pscircle[linewidth=0.06,dimen=outer](2.4,0.25){2.4}
\pscircle[linewidth=0.06,dimen=outer](5.43,1.88){1.25}
\pscircle[linewidth=0.06,dimen=outer](6.94,-0.27){1.5}
\pscircle[linewidth=0.06,dimen=outer](4.75,-1.8){1.33}
\end{pspicture} 
}
}
\qquad
\subfloat[\label{art3}]
{
\scalebox{1.2} 
{
\begin{pspicture*}(6.0,1.21)(4.24,-0.91)
\definecolor{color1789b}{rgb}{0.8,0.8,0.8}
\pscircle[linewidth=0.02,dimen=outer](2.4,0.25){2.4}
\pscircle[linewidth=0.02,dimen=outer](5.43,1.88){1.25}
\pscircle[linewidth=0.02,dimen=outer](6.94,-0.27){1.5}
\pscircle[linewidth=0.02,dimen=outer](4.75,-1.8){1.33}
\psframe[linewidth=0.02,dimen=middle](6.0,1.21)(4.24,-0.91)
\psline[linewidth=0.02cm,fillcolor=color1789b,linestyle=dashed,dash=0.1cm 0.05cm](4.7,0.87)(4.66,-0.49)
\psline[linewidth=0.02cm,fillcolor=color1789b,linestyle=dashed,dash=0.1cm 0.05cm](4.68,-0.49)(5.54,-0.75)
\psline[linewidth=0.02cm,fillcolor=color1789b,linestyle=dashed,dash=0.1cm 0.05cm](5.52,-0.73)(5.82,0.71)
\psline[linewidth=0.02cm,fillcolor=color1789b,linestyle=dashed,dash=0.1cm 0.05cm](5.82,0.71)(4.7,0.87)
\end{pspicture*} 
}

}
\caption[Four disks in a closed chain]
{
\label{fig:artistic}
{\bf Four disks in a closed chain.}  Suppose four disks have the contact graph shown in \subref{art1}.  An example of four such disks is shown in \subref{art2}.  If we zoom in on the hole formed by the four disks, we see the picture in \subref{art3}.  The angles inside of the dashed quadrilateral are clearly bigger than the corresponding overlap angles between the disks.  On the other hand, the sum of the angles inside of a quadrilateral is always exactly $360^\circ$.  So if four disks have the contact graph shown in \subref{art1}, then the sum of their overlap angles \emph{must} be less than $360^\circ$.
}
\end{figure}

\noindent Thus for instance there is no disk configuration having the following contact graph and overlap angles:

\begin{figure}[H]
\centering
\scalebox{1} 
{
\begin{pspicture}(0,-1.1829102)(5.2829103,1.1829102)
\definecolor{color1789b}{rgb}{0.8,0.8,0.8}
\psline[linewidth=0.02cm,fillcolor=color1789b,dotsize=0.07055555cm 2.0]{*-*}(1.0810156,0.5645117)(1.4810157,-0.83548826)
\psline[linewidth=0.02cm,fillcolor=color1789b,dotsize=0.07055555cm 2.0]{*-*}(1.4810157,-0.83548826)(3.8810155,-0.6354883)
\psline[linewidth=0.02cm,fillcolor=color1789b,dotsize=0.07055555cm 2.0]{*-*}(3.8810155,-0.6354883)(3.8810155,0.9645117)
\psline[linewidth=0.02cm,fillcolor=color1789b,dotsize=0.07055555cm 2.0]{*-*}(3.8810155,0.9645117)(1.0810156,0.5645117)
\usefont{T1}{ptm}{m}{n}
\rput(2.4724708,0.9895117){$100^\circ$}
\usefont{T1}{ptm}{m}{n}
\rput(4.392471,0.16951172){$100^\circ$}
\usefont{T1}{ptm}{m}{n}
\rput(2.7924707,-1.0104883){$100^\circ$}
\usefont{T1}{ptm}{m}{n}
\rput(0.8124707,-0.25048828){$100^\circ$}
\end{pspicture} 
}
\caption[An impossible contact graph given the indicated overlap angles]
{\label{fig:no here}
{\bf An impossible contact graph given the indicated overlap angles}, impossible because the sum of the indicated overlap angles is greater than $360^\circ$.
}
\end{figure}
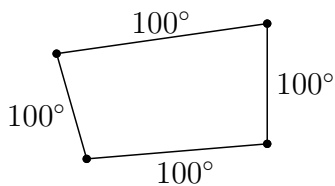

\noindent This might lead us to conclude that no disk configuration has a contact graph containing, as a sub-graph, the graph pictured in Figure \ref{fig:no here} including edge labels.  However, this conclusion does not turn out to be true.  Consider the following example:

\begin{figure}[H]
\centering
\subfloat[\label{fair1}]
{
\scalebox{1} 
{
\begin{pspicture}(0,-1.484)(4.3018947,1.484)
\definecolor{color1880b}{rgb}{0.8,0.8,0.8}
\pscircle[linestyle=none,fillstyle=solid,fillcolor=color1880b,linewidth=0.04,dimen=outer](2.18,0.38){1.07}
\pscircle[linestyle=none,fillstyle=solid,fillcolor=color1880b,linewidth=0.04,dimen=outer](2.26,-0.24){1.03}
\pscircle[linestyle=none,fillstyle=solid,fillcolor=color1880b,linewidth=0.04,dimen=outer](3.19,0.13){0.48}
\pscircle[linestyle=none,fillstyle=solid,fillcolor=color1880b,linewidth=0.04,dimen=outer](1.2,-0.16){0.65}
\pscircle[linewidth=0.04,dimen=outer](2.18,0.38){1.08}
\pscircle[linewidth=0.04,dimen=outer](2.26,-0.24){1.04}
\pscircle[linewidth=0.04,dimen=outer](3.19,0.13){0.49}
\pscircle[linewidth=0.04,dimen=outer](1.2,-0.16){0.66}
\psline[linewidth=0.0080cm,fillcolor=color1880b,arrowsize=0.05291667cm 2.0,arrowlength=1.4,arrowinset=0.4]{->}(1.12,0.5)(1.28,1.44)
\psline[linewidth=0.0080cm,fillcolor=color1880b,arrowsize=0.05291667cm 2.0,arrowlength=1.4,arrowinset=0.4]{->}(1.12,0.48)(0.0,0.36)
\psline[linewidth=0.0080cm,fillcolor=color1880b,arrowsize=0.05291667cm 2.0,arrowlength=1.4,arrowinset=0.4]{->}(1.38,-0.78)(0.38,-1.1)
\psline[linewidth=0.0080cm,fillcolor=color1880b,arrowsize=0.05291667cm 2.0,arrowlength=1.4,arrowinset=0.4]{->}(1.38,-0.78)(1.74,-1.48)
\psline[linewidth=0.0080cm,fillcolor=color1880b,arrowsize=0.05291667cm 2.0,arrowlength=1.4,arrowinset=0.4]{->}(3.28,-0.34)(3.26,-1.14)
\psline[linewidth=0.0080cm,fillcolor=color1880b,arrowsize=0.05291667cm 2.0,arrowlength=1.4,arrowinset=0.4]{->}(3.28,-0.34)(4.12,-0.18)
\psline[linewidth=0.0080cm,fillcolor=color1880b,arrowsize=0.05291667cm 2.0,arrowlength=1.4,arrowinset=0.4]{->}(3.22,0.62)(3.08,1.48)
\psline[linewidth=0.0080cm,fillcolor=color1880b,arrowsize=0.05291667cm 2.0,arrowlength=1.4,arrowinset=0.4]{->}(3.2,0.62)(4.06,0.56)
\psarc[linewidth=0.0080,fillcolor=color1880b](1.1,0.5){0.36}{75.96375}{191.30994}
\psarc[linewidth=0.0080,fillcolor=color1880b](3.2,0.6){0.36}{0.0}{96.340195}
\psarc[linewidth=0.0080,fillcolor=color1880b](1.38,-0.78){0.36}{196.69925}{296.56506}
\psarc[linewidth=0.0080,fillcolor=color1880b](3.26,-0.34){0.36}{270.0}{10.304847}
\usefont{T1}{ptm}{m}{n}
\rput(0.6814551,0.965){$100^\circ$}
\usefont{T1}{ptm}{m}{n}
\rput(1.141455,-1.365){$100^\circ$}
\usefont{T1}{ptm}{m}{n}
\rput(3.8214552,1.045){$100^\circ$}
\usefont{T1}{ptm}{m}{n}
\rput(3.8814551,-0.715){$100^\circ$}
\end{pspicture} 
}
}\qquad
\subfloat[\label{fair2}]
{
\scalebox{1} 
{
\begin{pspicture}(0,-1.0566211)(4.6450157,1.0891992)
\definecolor{color2225b}{rgb}{0.8,0.8,0.8}
\psline[linewidth=0.0080cm,fillcolor=color2225b,dotsize=0.07055555cm 2.0]{*-*}(2.0410156,0.9508008)(0.041015625,-0.24919921)
\psline[linewidth=0.0080cm,fillcolor=color2225b,dotsize=0.07055555cm 2.0]{*-*}(0.041015625,-0.24919921)(2.4410157,-1.0491992)
\psline[linewidth=0.0080cm,fillcolor=color2225b,dotsize=0.07055555cm 2.0]{*-*}(2.4410157,-1.0491992)(4.6410155,0.35080078)
\psline[linewidth=0.0080cm,fillcolor=color2225b,dotsize=0.07055555cm 2.0]{*-*}(4.6410155,0.35080078)(2.0410156,0.9508008)
\psline[linewidth=0.0080cm,fillcolor=color2225b,dotsize=0.07055555cm 2.0]{*-*}(2.0410156,0.9508008)(2.4410157,-1.0491992)
\usefont{T1}{ptm}{m}{n}
\rput(2.6524706,0.03580078){$130^\circ$}
\usefont{T1}{ptm}{m}{n}
\rput(0.9524707,-0.8841992){$100^\circ$}
\usefont{T1}{ptm}{m}{n}
\rput(3.7524707,-0.60419923){$100^\circ$}
\usefont{T1}{ptm}{m}{n}
\rput(3.5324707,0.89580077){$100^\circ$}
\usefont{T1}{ptm}{m}{n}
\rput(0.8124707,0.5758008){$100^\circ$}
\end{pspicture} 
}
}
\caption[Another disk configuration, its contact graph, and overlap angles]
{
\label{fig:fair}
{\bf Another disk configuration, its contact graph, and overlap angles.}
}
\end{figure}
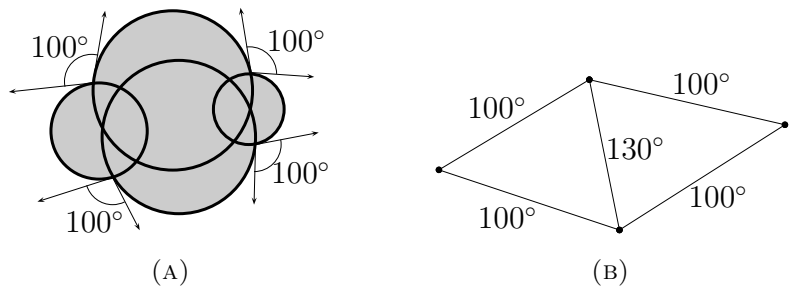%
At this point, we give up on the question of existence and summarize our findings on this question in the following observation:

\begin{observation*}
Given a graph $G$ with edges labeled with angles, it is hard to determine whether there is any disk configuration which has $G$ as its contact graph, with those overlap angles.
\end{observation*}

\section{Rigidity of thin disk configurations}

Even though we can't easily figure out whether a given pattern (a contact graph, having edges labeled with angles) can actually be realized by a disk configuration, we can still ask what properties a disk configuration should have to guarantee that it is the \emph{only one} having its particular pattern:

\begin{question*}
Suppose we start with a disk configuration.  What conditions guarantee that there is essentially no other disk configuration having the same contact graph and overlap angles?
\end{question*}

\noindent What it means for two disk configurations to be \emph{similar} is the same as what it meant for two circle packings to be similar.  Then, as before, the following is not too hard to believe:

\begin{fact*}
Suppose that the disk configuration $\calC$ consists of only finitely many disks.  Then there are disk configurations having the same ``pattern'' as $\calC$, but which are not similar to $\calC$.
\end{fact*}

\begin{exercise}
\label{simdc}
Draw a disk configuration which has the same contact graph and overlap angles as the configuration depicted in Figure \ref{jeff}, but which is not similar to it.
\end{exercise}

\noindent Also, as with circle packings, there are infinite disk configurations which have the same contact graph and overlap angles, but which are not similar:

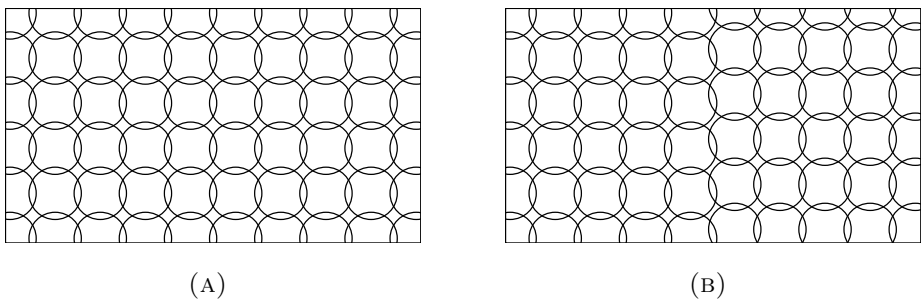
\begin{figure}[H]
\centering
\subfloat[\label{sqq1}]
{
\scalebox{.6}
{
\begin{pspicture*}(-.1,-.1)(9.1,5.1)
\psframe[dimen=middle](-.1,-.1)(9.1,5.1)
\pscircle[dimen=middle](0,0){.58}
\pscircle[dimen=middle](1,0){.58}
\pscircle[dimen=middle](2,0){.58}
\pscircle[dimen=middle](3,0){.58}
\pscircle[dimen=middle](4,0){.58}
\pscircle[dimen=middle](5,0){.58}
\pscircle[dimen=middle](6,0){.58}
\pscircle[dimen=middle](7,0){.58}
\pscircle[dimen=middle](8,0){.58}
\pscircle[dimen=middle](9,0){.58}%

\pscircle[dimen=middle](0,1){.58}
\pscircle[dimen=middle](1,1){.58}
\pscircle[dimen=middle](2,1){.58}
\pscircle[dimen=middle](3,1){.58}
\pscircle[dimen=middle](4,1){.58}
\pscircle[dimen=middle](5,1){.58}
\pscircle[dimen=middle](6,1){.58}
\pscircle[dimen=middle](7,1){.58}
\pscircle[dimen=middle](8,1){.58}
\pscircle[dimen=middle](9,1){.58}%

\pscircle[dimen=middle](0,2){.58}
\pscircle[dimen=middle](1,2){.58}
\pscircle[dimen=middle](2,2){.58}
\pscircle[dimen=middle](3,2){.58}
\pscircle[dimen=middle](4,2){.58}
\pscircle[dimen=middle](5,2){.58}
\pscircle[dimen=middle](6,2){.58}
\pscircle[dimen=middle](7,2){.58}
\pscircle[dimen=middle](8,2){.58}
\pscircle[dimen=middle](9,2){.58}%

\pscircle[dimen=middle](0,3){.58}
\pscircle[dimen=middle](1,3){.58}
\pscircle[dimen=middle](2,3){.58}
\pscircle[dimen=middle](3,3){.58}
\pscircle[dimen=middle](4,3){.58}
\pscircle[dimen=middle](5,3){.58}
\pscircle[dimen=middle](6,3){.58}
\pscircle[dimen=middle](7,3){.58}
\pscircle[dimen=middle](8,3){.58}
\pscircle[dimen=middle](9,3){.58}%

\pscircle[dimen=middle](0,4){.58}
\pscircle[dimen=middle](1,4){.58}
\pscircle[dimen=middle](2,4){.58}
\pscircle[dimen=middle](3,4){.58}
\pscircle[dimen=middle](4,4){.58}
\pscircle[dimen=middle](5,4){.58}
\pscircle[dimen=middle](6,4){.58}
\pscircle[dimen=middle](7,4){.58}
\pscircle[dimen=middle](8,4){.58}
\pscircle[dimen=middle](9,4){.58}%

\pscircle[dimen=middle](0,5){.58}
\pscircle[dimen=middle](1,5){.58}
\pscircle[dimen=middle](2,5){.58}
\pscircle[dimen=middle](3,5){.58}
\pscircle[dimen=middle](4,5){.58}
\pscircle[dimen=middle](5,5){.58}
\pscircle[dimen=middle](6,5){.58}
\pscircle[dimen=middle](7,5){.58}
\pscircle[dimen=middle](8,5){.58}
\pscircle[dimen=middle](9,5){.58}
\end{pspicture*}
}
}
\qquad
\subfloat[\label{sqq2}]
{
\scalebox{.6}
{
\begin{pspicture*}(-.1,-.1)(9.1,5.1)
\psframe[dimen=middle](-.1,-.1)(9.1,5.1)
\pscircle[dimen=middle](0,0){.58}
\pscircle[dimen=middle](1,0){.58}
\pscircle[dimen=middle](2,0){.58}
\pscircle[dimen=middle](3,0){.58}
\pscircle[dimen=middle](4,0){.58}
\pscircle[dimen=middle](4.978,.2){.58}
\pscircle[dimen=middle](5.978,.2){.58}
\pscircle[dimen=middle](6.978,.2){.58}
\pscircle[dimen=middle](7.978,.2){.58}
\pscircle[dimen=middle](8.978,.2){.58}%

\pscircle[dimen=middle](0,1){.58}
\pscircle[dimen=middle](1,1){.58}
\pscircle[dimen=middle](2,1){.58}
\pscircle[dimen=middle](3,1){.58}
\pscircle[dimen=middle](4,1){.58}
\pscircle[dimen=middle](4.978,1.2){.58}
\pscircle[dimen=middle](5.978,1.2){.58}
\pscircle[dimen=middle](6.978,1.2){.58}
\pscircle[dimen=middle](7.978,1.2){.58}
\pscircle[dimen=middle](8.978,1.2){.58}%

\pscircle[dimen=middle](0,2){.58}
\pscircle[dimen=middle](1,2){.58}
\pscircle[dimen=middle](2,2){.58}
\pscircle[dimen=middle](3,2){.58}
\pscircle[dimen=middle](4,2){.58}
\pscircle[dimen=middle](4.978,2.2){.58}
\pscircle[dimen=middle](5.978,2.2){.58}
\pscircle[dimen=middle](6.978,2.2){.58}
\pscircle[dimen=middle](7.978,2.2){.58}
\pscircle[dimen=middle](8.978,2.2){.58}%

\pscircle[dimen=middle](0,3){.58}
\pscircle[dimen=middle](1,3){.58}
\pscircle[dimen=middle](2,3){.58}
\pscircle[dimen=middle](3,3){.58}
\pscircle[dimen=middle](4,3){.58}
\pscircle[dimen=middle](4.978,3.2){.58}
\pscircle[dimen=middle](5.978,3.2){.58}
\pscircle[dimen=middle](6.978,3.2){.58}
\pscircle[dimen=middle](7.978,3.2){.58}
\pscircle[dimen=middle](8.978,3.2){.58}%

\pscircle[dimen=middle](0,4){.58}
\pscircle[dimen=middle](1,4){.58}
\pscircle[dimen=middle](2,4){.58}
\pscircle[dimen=middle](3,4){.58}
\pscircle[dimen=middle](4,4){.58}
\pscircle[dimen=middle](4.978,4.2){.58}
\pscircle[dimen=middle](5.978,4.2){.58}
\pscircle[dimen=middle](6.978,4.2){.58}
\pscircle[dimen=middle](7.978,4.2){.58}
\pscircle[dimen=middle](8.978,4.2){.58}%

\pscircle[dimen=middle](0,5){.58}
\pscircle[dimen=middle](1,5){.58}
\pscircle[dimen=middle](2,5){.58}
\pscircle[dimen=middle](3,5){.58}
\pscircle[dimen=middle](4,5){.58}
\pscircle[dimen=middle](4.978,5.2){.58}
\pscircle[dimen=middle](5.978,5.2){.58}
\pscircle[dimen=middle](6.978,5.2){.58}
\pscircle[dimen=middle](7.978,5.2){.58}
\pscircle[dimen=middle](8.978,5.2){.58}
\end{pspicture*}
}
}
\caption[Two infinite disk configurations which are not similar, but which have the same contact graph and overlap angles]
{
\label{fig:sqq}
{\bf Two infinite disk configurations which are not similar, but which have the same contact graph and overlap angles.}
}
\end{figure}

\noindent Similarly to before, it will help to consider \emph{triangulated} disk configurations: a disk configuration is \emph{triangulated} if its contact graph ``cuts the plane into triangles,'' in the same way as for circle packings.  Then the following theorem turns out to be true:

\begin{maint}
If two thin triangulated disk configurations have the same pattern (contact graph, including overlap angles), then the two disk configurations are similar.
\end{maint}

\noindent We still need to define the extra adjective \emph{thin} in the statement of the above theorem:

\begin{definition*}
A disk configuration is \bfit{thin} if no three of its disks meet at any point.
\end{definition*}

\begin{figure}[H]
\centering
\subfloat[\label{yes thin}]
{
\scalebox{.6} 
{
\begin{pspicture}(0,-1.69)(4.16,1.69)
\pscircle[linewidth=0.04,dimen=outer](1.32,0.21){1.32}
\pscircle[linewidth=0.04,dimen=outer](3.23,0.76){0.93}
\pscircle[linewidth=0.04,dimen=outer](3.01,-0.84){0.85}
\end{pspicture} 
}
}
\qquad
\subfloat[\label{not thin}]
{
\scalebox{.6} 
{
\begin{pspicture}(0,-1.51)(3.84,1.51)
\pscircle[linewidth=0.04,dimen=outer](1.32,0.19){1.32}
\pscircle[linewidth=0.04,dimen=outer](2.91,0.44){0.93}
\pscircle[linewidth=0.04,dimen=outer](2.57,-0.66){0.85}
\end{pspicture} 
}
}
\caption[A disk configuration which is thin and one which is not]
{
{\bf A disk configuration which is thin and one which is not.}  The configuration in \subref{not thin} is not thin, because there are points which are in all three of the disks at once.  This does not happen in \subref{yes thin}, so the configuration in \subref{yes thin} is thin.
}
\end{figure}
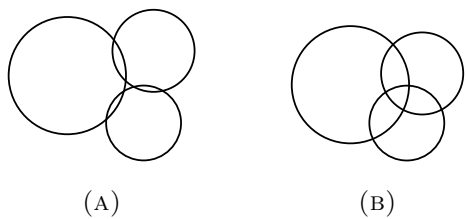

\noindent We give two examples of disk configurations to which the above uniqueness theorem applies:

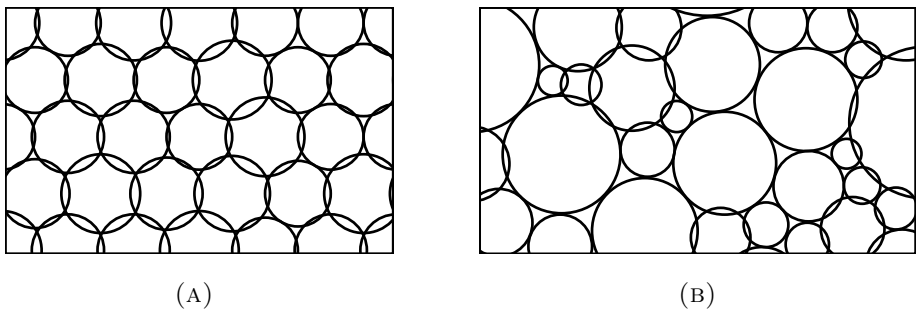
\begin{figure}[H]
\centering
\subfloat[]
{
\scalebox{0.436797}{
\begin{pspicture*}(0.1,-0.1)(11.9,7.4)
\psframe[dimen=middle,linewidth=0.0824](0.1,-0.1)(11.9,7.4)
\pscircle[dimen=middle,linewidth=0.0824](00,0){1.2}
\pscircle[dimen=middle,linewidth=0.0824](02,0){1.1}
\pscircle[dimen=middle,linewidth=0.0824](04,0){1.1}
\pscircle[dimen=middle,linewidth=0.0824](06,0){1.2}
\pscircle[dimen=middle,linewidth=0.0824](08,0){1}
\pscircle[dimen=middle,linewidth=0.0824](10,0){1.1}
\pscircle[dimen=middle,linewidth=0.0824](12,0){1.1}%

\pscircle[dimen=middle,linewidth=0.0824](01,1.73){1.05}
\pscircle[dimen=middle,linewidth=0.0824](03,1.73){1.2}
\pscircle[dimen=middle,linewidth=0.0824](05,1.73){1.1}
\pscircle[dimen=middle,linewidth=0.0824](07,1.73){1.15}
\pscircle[dimen=middle,linewidth=0.0824](09,1.73){1}
\pscircle[dimen=middle,linewidth=0.0824](11,1.73){1.2}
\pscircle[dimen=middle,linewidth=0.0824](13,1.73){1.2}%

\pscircle[dimen=middle,linewidth=0.0824](00,3.46){1}
\pscircle[dimen=middle,linewidth=0.0824](02,3.46){1.1}
\pscircle[dimen=middle,linewidth=0.0824](04,3.46){1.1}
\pscircle[dimen=middle,linewidth=0.0824](06,3.46){1}
\pscircle[dimen=middle,linewidth=0.0824](08,3.46){1.2}
\pscircle[dimen=middle,linewidth=0.0824](10,3.46){1}
\pscircle[dimen=middle,linewidth=0.0824](12,3.46){1}%

\pscircle[dimen=middle,linewidth=0.0824](01,5.19){1}
\pscircle[dimen=middle,linewidth=0.0824](03,5.19){1.1}
\pscircle[dimen=middle,linewidth=0.0824](05,5.19){1}
\pscircle[dimen=middle,linewidth=0.0824](07,5.19){1.2}
\pscircle[dimen=middle,linewidth=0.0824](09,5.19){1}
\pscircle[dimen=middle,linewidth=0.0824](11,5.19){1.1}
\pscircle[dimen=middle,linewidth=0.0824](13,5.19){1.1}%

\pscircle[dimen=middle,linewidth=0.0824](00,6.928){1.13}
\pscircle[dimen=middle,linewidth=0.0824](02,6.928){1}
\pscircle[dimen=middle,linewidth=0.0824](04,6.928){1.14}
\pscircle[dimen=middle,linewidth=0.0824](06,6.928){1.15}
\pscircle[dimen=middle,linewidth=0.0824](08,6.928){1}
\pscircle[dimen=middle,linewidth=0.0824](10,6.928){1}
\pscircle[dimen=middle,linewidth=0.0824](12,6.928){1}%
\end{pspicture*}
}
}
\qquad
\subfloat[]
{
\scalebox{.9} 
{
\begin{pspicture*}(7.56,1.12)(1.12,-2.52)
\psframe[linewidth=0.04,dimen=middle](7.56,1.12)(1.12,-2.52)
\pscircle[linewidth=0.04,dimen=outer](1.01,0.27){1.01}
\pscircle[linewidth=0.04,dimen=outer](2.57,0.83){0.67}
\pscircle[linewidth=0.04,dimen=outer](3.37,-0.07){0.65}
\pscircle[linewidth=0.04,dimen=outer](2.2,0.04){0.24}
\pscircle[linewidth=0.04,dimen=outer](2.62,-0.02){0.32}
\pscircle[linewidth=0.04,dimen=outer](2.33,-1.05){0.89}
\pscircle[linewidth=0.04,dimen=outer](3.62,0.82){0.5}
\pscircle[linewidth=0.04,dimen=outer](4.56,0.28){0.72}
\pscircle[linewidth=0.04,dimen=outer](3.6,-0.98){0.42}
\pscircle[linewidth=0.04,dimen=outer](4.03,-0.49){0.25}
\pscircle[linewidth=0.04,dimen=outer](1.01,-1.19){0.57}
\pscircle[linewidth=0.04,dimen=outer](4.74,-1.18){0.78}
\pscircle[linewidth=0.04,dimen=outer](3.56,-2.18){0.8}
\pscircle[linewidth=0.04,dimen=outer](1.4,-2.06){0.52}
\pscircle[linewidth=0.04,dimen=outer](2.31,-2.41){0.49}
\pscircle[linewidth=0.04,dimen=outer](4.68,-2.28){0.46}
\pscircle[linewidth=0.04,dimen=outer](5.35,-2.09){0.35}
\pscircle[linewidth=0.04,dimen=outer](5.98,-1.52){0.54}
\pscircle[linewidth=0.04,dimen=outer](5.94,-0.24){0.78}
\pscircle[linewidth=0.04,dimen=outer](6.33,0.83){0.39}
\pscircle[linewidth=0.04,dimen=outer](6.79,0.35){0.29}
\pscircle[linewidth=0.04,dimen=outer](7.7,-0.58){1.14}
\pscircle[linewidth=0.04,dimen=outer](6.54,-1.04){0.24}
\pscircle[linewidth=0.04,dimen=outer](6.77,-1.51){0.29}
\pscircle[linewidth=0.04,dimen=outer](6.63,-2.15){0.49}
\pscircle[linewidth=0.04,dimen=outer](5.96,-2.38){0.34}
\pscircle[linewidth=0.04,dimen=outer](7.27,-1.85){0.33}
\pscircle[linewidth=0.04,dimen=outer](7.42,1.2){0.88}
\pscircle[linewidth=0.04,dimen=outer](5.53,0.89){0.45}
\pscircle[linewidth=0.04,dimen=outer](7.36,-2.7){0.56}
\pscircle[linewidth=0.04,dimen=outer](5.29,-3.01){0.61}
\pscircle[linewidth=0.04,dimen=outer](4.48,2.3){1.32}
\end{pspicture*} 
}
}
\caption[Pieces of thin triangulated disk configurations]
{
{\bf Pieces of thin triangulated disk configurations.}
}
\end{figure}

\noindent We also give an example of a disk configuration which is not covered by the above theorem:

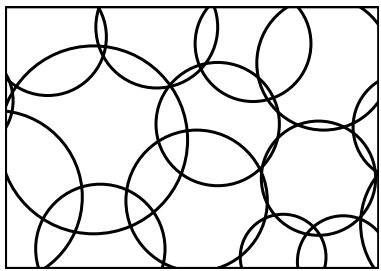
\begin{figure}[H]
\centering
\scalebox{1} 
{
\begin{pspicture*}(6.3,1.77)(1.34,-1.71)
\pscircle[linewidth=0.04,dimen=outer](1.21,-0.76){1.17}
\pscircle[linewidth=0.04,dimen=outer](2.51,0.0){1.27}
\pscircle[linewidth=0.04,dimen=outer](2.6,-1.45){0.88}
\pscircle[linewidth=0.04,dimen=outer](3.88,-0.81){0.96}
\pscircle[linewidth=0.04,dimen=outer](4.16,0.21){0.84}
\pscircle[linewidth=0.04,dimen=outer](5.5,-0.51){0.78}
\pscircle[linewidth=0.04,dimen=outer](5.57,1.02){0.91}
\pscircle[linewidth=0.04,dimen=outer](3.35,1.5){0.83}
\pscircle[linewidth=0.04,dimen=outer](4.63,1.28){0.79}
\pscircle[linewidth=0.04,dimen=outer](1.9,1.37){0.8}
\pscircle[linewidth=0.04,dimen=outer](0.73,0.52){0.73}
\pscircle[linewidth=0.04,dimen=outer](5.03,-1.56){0.59}
\pscircle[linewidth=0.04,dimen=outer](6.68,0.17){0.74}
\pscircle[linewidth=0.04,dimen=outer](7.11,-1.12){1.07}
\pscircle[linewidth=0.04,dimen=outer](5.83,-1.62){0.63}
\psframe[linewidth=0.04,dimen=middle](6.3,1.77)(1.34,-1.71)
\end{pspicture*}
}
\caption[A triangulated, but non-thin, disk configuration]
{
{\bf A triangulated, but non-thin, disk configuration.}
}
\end{figure}

Another way of stating the uniqueness theorem above is to say that thin triangulated disk configurations are \emph{rigid}, meaning that we cannot modify them while keeping the same contact graph and overlap angles, except by applying rotations, scaling, etc.

As before, consideration of Figure \ref{fig:sqq} makes the following believable:

\begin{fact*}
Suppose that $\calC$ is a thin disk configuration which is not triangulated.  Then there are disk configurations having the same pattern as $\calC$, but which are not similar to $\calC$.
\end{fact*}

\noindent Thus the uniqueness guaranteed by the above theorem is the best uniqueness we could hope for, at least for thin disk configurations.

I proved the uniqueness theorem for triangulated thin disk configurations as part of my Ph.\ D.\ thesis \cite{mishchenko-thesis}.  The ``thinness'' condition is used strongly in the proof, but it seems likely that it is not necessary, so we can speculate:

\begin{conjecture*}
The uniqueness theorem for thin triangulated disk configurations is still true if we completely get rid of the thinness requirement from the statement.
\end{conjecture*}

\noindent A \emph{conjecture} is a mathematical statement that is thought by some to be true, but which has not yet been proved.  In 1999 Z.-X.\ He, in \cite{MR1680531}, proved the following:

\begin{theorem*}
The above uniqueness theorem is still true if we get rid of the thinness requirement from the statement, provided that none of the overlap angles is bigger than $90^\circ$.
\end{theorem*}

\noindent However, the restriction that overlap angles stay below $90^\circ$ is a strong one, and He never published a version of this theorem where he managed to eliminate this restriction.

\section{Closing remarks and further references}

You might wonder what real-world applications circle packing has.  To describe examples in any detail would require considerably more room than we have here.  However, for instance, some people are trying to use circle packings as a tool in medical imaging, for example to map regions of the brain, c.f.\ \cite{MR2131318}*{Section 23.4}.

Stephenson has written a book \cite{MR2131318} which gives an introduction to circle packing, including many drawings and examples.  It should be readable by an advanced undergraduate math major.  For more advanced readers, a nice survey of the research area as of 2011 is contained in the first half of \cite{MR2884870}, with a focus on the contributions of Oded Schramm.  There is also a survey of the area including many references in my thesis \cite{mishchenko-thesis}*{Chapter 2}.


\end{document}